\documentclass[11pt]{article}
\usepackage[a4paper,margin=27mm]{geometry}
\usepackage{amsmath,amssymb,amsthm,mathtools}
\usepackage{xcolor,graphicx,booktabs}
\usepackage[hidelinks]{hyperref}
\newcommand{\N}{\mathbb N}
\newcommand{\Z}{\mathbb Z}
\newcommand{\F}{\mathbb F}
\newcommand{\X}{\boldsymbol X}
\newcommand{\x}{\boldsymbol x}
\newcommand{\vct}{\boldsymbol v}
\newcommand{\rct}{\boldsymbol r}
\newcommand{\sct}{\boldsymbol s}
\newcommand{\one}{\mathbf1}
\newcommand{\kerp}{\mathcal K_p}
\DeclareMathOperator{\adj}{adj}
\newtheorem{theorem}{Theorem}[section]
\newtheorem{lemma}[theorem]{Lemma}
\newtheorem{proposition}[theorem]{Proposition}
\newtheorem{corollary}[theorem]{Corollary}
\theoremstyle{definition}
\newtheorem{definition}[theorem]{Definition}
\newtheorem*{definition*}{Definition}
\theoremstyle{remark}
\newtheorem{remark}[theorem]{Remark}
\graphicspath{{images/}}
\definecolor{zero}{RGB}{247,245,239}
\definecolor{one}{RGB}{40,116,142}
\definecolor{two}{RGB}{204,91,70}
\newcommand{\legendthree}{\colorbox{zero}{\strut\phantom{M}} $0$\quad
  \colorbox{one}{\strut\phantom{M}} $1$\quad
  \colorbox{two}{\strut\phantom{M}} $2$}
\title{Automatic $n$-dimensional recurrences\\
over finite $p$-groups\\
with negative predecessor coordinates}
\author{Bogdan Dumitru\footnote{Faculty of Mathematics and Computer Science, University of Bucharest, Academiei 14, Bucharest (RO-010014), Romania.}\,\,,\,\,Mihai Prunescu\footnote{Research Center for Logic, Optimization and
Security (LOS), Faculty of Mathematics and Computer Science, University of
Bucharest, Academiei 14, Bucharest (RO-010014), Romania.}\,\,\footnote{Simion
Stoilow Institute of Mathematics of the Romanian Academy, Research unit 5,
P. O. Box 1-764, Bucharest (RO-014700), Romania.
E-mail: \texttt{mihai.prunescu@imar.ro}.}}
\date{}
\begin{document}
\maketitle

\begin{abstract}
We study $n$-dimensional recurrences on $\N^n$ with values in a finite
$p$-group. An array is $p$-automatic if a finite automaton computes its
entries from the base-$p$ digits of all $n$ coordinates. The recurrence
rule is an ordered product of group endomorphisms, with a prescribed
$p$-automatic additional sequence and compatible $p$-automatic boundary
data. The predecessor set is finite, nonempty and lexicographically
positive, meaning that each vector's first nonzero coordinate is positive.
Under these assumptions, the automaticity question depends on
the dimension. For $n=1$, lexicographically positive predecessor
offsets are positive integers, and every such recurrence is
$p$-automatic. For $n=2$, every lexicographically positive predecessor
set also gives a $p$-automatic array, including sets with negative
second coordinates. For $n\ge3$, lexicographic positivity alone is
insufficient. We prove automaticity when at most one coordinate of
the predecessor set contains entries of both signs. Other coordinates
may contain negative entries if they contain no positive entries.
For nonnegative predecessors, we give a direct automaton construction
and an explicit bound on its number of states in the abelian case.
The proof allowing negative entries uses algebraic power-series
division, lifting through powers of $p$, central subgroups, and integer
changes of coordinates that preserve the prescribed boundary.
A fixed three-dimensional recurrence with one nonzero boundary value
fails to be $p$-automatic for every prime $p$. For a two-parameter family
with a specified rule and a single boundary seed, we give an exact
automaticity criterion modulo prime powers.
Nonabelian examples include figures computed from their recurrences.
Certificates consisting of transition tables and algebraic identities
verify further automata and prove a previously proposed formula
counting nonzero Delannoy residues modulo $27$.
\end{abstract}

\section{Introduction and statements}\label{sec:introduction}

An $n$-dimensional recurrence determines an array indexed by $\N^n$
from prescribed boundary values and values at predecessor points.
We ask when the resulting array is $p$-automatic, meaning that a
deterministic automaton with finitely many states computes an entry by
reading the base-$p$ digits of its $n$ coordinates together. The
existence of a recurrence with finitely many predecessor offsets and
automatic boundary data does not by itself imply automaticity of the
solution. Both the algebraic rule and the geometry of the predecessor
offsets matter.

Prunescu \cite[Theorem~5.2]{Prunescu} proves automaticity for
$\F_p$-affine recurrences over finite fields $\F_q$ of characteristic
$p$ when predecessor vectors are coordinatewise nonnegative and boundary
faces are automatic. He asks \cite[Section~7, Open problems~3--5]{Prunescu}
whether this remains true for predecessor sets with negative entries,
for homomorphisms of arbitrary finite abelian $p$-groups, and for
both together. We allow finite $p$-groups, meaning finite groups of order a power of
$p$, and predecessor vectors with negative entries in one common
coordinate. Theorem~\ref{thm:mixed} then permits negative
entries in further coordinates, provided those coordinates have no
positive entries. The rule can be an ordered product of endomorphisms,
where an endomorphism is a homomorphism from the group to itself.
Theorems~\ref{thm:abelian} and~\ref{thm:groups} state the resulting
automaticity conclusions.

For these group rules, with a $p$-automatic additional sequence and
compatible $p$-automatic boundary data, the dimensions separate as follows.
In dimension one, lexicographic positivity requires positive integer
offsets. The boundary is a finite initial segment, and the nonnegative
construction gives automaticity. Section~\ref{sec:one-dimension}
states this case and explains the stronger conclusion of eventual
periodicity when the additional sequence is constant. In dimension
two, the first coordinate of every lexicographically positive vector
is nonnegative. Negative entries can therefore occur only in the
second coordinate, so Corollary~\ref{cor:plane} covers every such
predecessor set. In dimensions at least three, negative entries can
occur in several coordinates that also contain positive entries.
The counterexample in Section~\ref{sec:limit} shows that lexicographic
positivity then ceases to be a sufficient condition. Theorem~\ref{thm:mixed}
still gives automaticity when at most one coordinate contains both
signs. In this sense, three is the first possible dimension of failure
under the stated algebraic and data assumptions.

For Prunescu's questions, Theorem~\ref{thm:abelian} gives a positive
answer to Open problem~4, retaining the nonnegative predecessor
hypothesis of his theorem. Corollary~\ref{cor:plane} answers Open
problems~3 and~5 affirmatively in dimension two, and
Theorem~\ref{thm:mixed} does so in every dimension under its sign
condition. Section~\ref{sec:limit} gives counterexamples to the
unrestricted assertions in Open problems~3 and~5 in every dimension
$n\ge3$, already for $\F_p$-linear rules over $\F_p$.

Theorem~\ref{thm:time-boundary} treats another sufficient family with
negative entries in several coordinates. Its boundary consists of
the $d$ coordinate slices $\{x_1=j\}$ for $0\le j<d$. Each slice
is a copy of $\N^{n-1}$ and is infinite when $n\ge2$. An integer
change of coordinates reduces this recurrence to the nonnegative case.

The proof treats the boundary differently according to the predecessor
geometry. For nonnegative vectors, the boundary contribution is known
and the solution is an automatic series divided by a polynomial with
constant coefficient $1$. We give its automaton by extending the
coefficient selection method of Rowland--Yassawi \cite{RY15} from
polynomial to automatic numerators. With negative entries,
the generating-function equation also contains unknown boundary series.
When the negative entries are confined to one coordinate, algebraic
Weierstrass division determines those series over a finite field.
Reduction modulo successive powers of $p$ gives the abelian group
case. A central subgroup then reduces a nonabelian group recurrence
to a smaller group and an abelian recurrence with an automatic affine term.

Theorem~\ref{thm:verified-search} gives a second route to an explicit
automaton. It samples arrays obtained by fixing the low digits of the
indices, called kernel arrays, to construct a candidate. It then checks
the recurrence and boundary values at all reachable states. Automaticity
guarantees eventual termination, while each accepted candidate has a
certificate valid at all indices. A matrix that counts automaton
transitions also proves Laurberg's proposed formula
\cite[Eq.~(10.47)]{Laurberg} for the number of nonzero Delannoy residues
modulo $27$ in squares of side $3^r$, for every $r\ge4$.

The scalar algebraic
principle for one possible negative direction appears in
\cite[Theorem~13]{BMP}, stated over $\mathbb C$, and is interpreted
through division in \cite{HK}. Both papers work in characteristic zero.
We use division valid in arbitrary
characteristic \cite{Rond} and supply the matrix reduction needed for
coupled recurrences. The group induction uses the decomposition method
of \cite[Theorem~4 and Proposition~8(3)]{Pivato}. Related
cellular-automaton results appear in \cite{GNW,RYCA,Nesme}.
Nesme \cite[Lemma~4]{Nesme} also proves automaticity of convolutions
taking values in a finite commutative monoid when each coefficient
receives only finitely many nonzero contributions. We need
only multiplication by polynomials, which reduces to finite sums of
translations by Lemma~\ref{lem:closure}.
The obstruction adapts the recurrence in
\cite[Section~4.5, Example~6]{BMP}.

Fix a prime $p$ and put $\N=\{0,1,\ldots\}$. Coordinatewise
inequalities are used unless lexicographic order is specified.
A vector is \emph{lexicographically positive} if its first nonzero
coordinate is positive. Let
\[
 P=\{\vct_1,\ldots,\vct_m\}\subset\Z^n,\qquad
 \vct_j>_{\rm lex}0,\qquad d_i=\max(0,\max_j(v_j)_i),
\]
where $P$ is finite and nonempty. A predecessor of $\x$ is a point
$\x-\vct_j$. The interior and prescribed boundary are
\[
 I=\boldsymbol d+\N^n,\qquad
C=\N^n\setminus I
   =\bigcup_{i=1}^n\bigcup_{0\le j<d_i}\{x_i=j\}.
\]
At an interior point all required predecessors belong to $\N^n$.
At a boundary point at least one required predecessor leaves $\N^n$.
An empty range $0\le j<d_i$ contributes no face.
There are $\sum_i d_i$ faces. Each face has infinitely many
points when $n\ge2$, while for $n=1$ the boundary consists of $d_1$
points.

\begin{definition}\label{def:automatic}
A function on $\N^n$ with a finite output alphabet is
\emph{$p$-automatic} if a finite deterministic automaton with output
computes its value by reading the base-$p$ digits of all coordinates
together, padded to a common length. The output must be invariant under
additional zero padding. A boundary function $c:C\to A$ is automatic
on its faces if its restriction to each face $x_i=j$, $0\le j<d_i$,
is an $(n-1)$-dimensional $p$-automatic function. The values on
intersecting faces must agree. A face in dimension one is a single point.
\end{definition}

The additional sequence $b$ is a prescribed array on $\N^n$.
In additive notation, $b(\x)$ is the affine term of the recurrence.
In a rule written multiplicatively, it is the prescribed factor on
the right of the ordered product.

Our sufficient geometric condition is
\begin{equation}\label{eq:geometry}
 \text{there exists }k\text{ such that }(v_j)_i\ge0
 \text{ for all }j\text{ and all }i\ne k.
\end{equation}
The coordinate $k$ is common to the whole predecessor set. For example,
$\{(1,-2),(0,1)\}$ satisfies the condition with $k=2$.
The set $\{(1,2,-1),(1,-1,2),(2,1,1)\}$ has negative entries
in two different coordinates and does not satisfy the condition.
Throughout the statements below, the predecessor set is also assumed
to be finite, nonempty and lexicographically positive as specified above.

\begin{theorem}[Abelian case]\label{thm:abelian}
Let $H$ be a finite abelian $p$-group, $h:H^m\to H$ a homomorphism,
and $b:\N^n\to H$ a $p$-automatic function. Suppose $P$ satisfies
\eqref{eq:geometry} and $c:C\to H$ is automatic on its faces with
compatible values on intersections. Then the unique solution of
\[
 a(\x)=c(\x)\quad(\x\in C),\qquad
 a(\x)=h(a(\x-\vct_1),\ldots,a(\x-\vct_m))+b(\x)
 \quad(\x\in I)
\]
is $p$-automatic. For coordinatewise nonnegative predecessors, the proof
gives an explicit terminating automaton construction from the rule and
automata for the data. The term $b$ may in particular be constant.
\end{theorem}

\begin{theorem}[Finite $p$-groups]\label{thm:groups}
Let $G$ be a finite $p$-group and let $\alpha_1,\ldots,\alpha_m$
be endomorphisms of $G$. Suppose $P$ satisfies \eqref{eq:geometry},
$b:\N^n\to G$ is $p$-automatic, and $c:C\to G$ is automatic on
its faces with compatible values on intersections. Then the unique
solution of
\[
 a(\x)=c(\x)\quad(\x\in C),\qquad
 a(\x)=\alpha_1(a(\x-\vct_1))\cdots
         \alpha_m(a(\x-\vct_m))b(\x)\quad(\x\in I)
\]
is $p$-automatic. The order of the factors is fixed as displayed.
\end{theorem}

\begin{corollary}[Direct-product homomorphism rules]\label{cor:hom}
Under the same geometric and data assumptions, the conclusion holds
for the rule $a(\x)=h(a(\x-\vct_1),\ldots,a(\x-\vct_m))b(\x)$
when $h:G^m\to G$ is a group homomorphism.
\end{corollary}
\begin{proof}
Define $\alpha_j(g)=h(1,\ldots,g,\ldots,1)$ using the $j$th input.
Each $\alpha_j$ is an endomorphism, and the product decomposition
of a tuple in $G^m$ gives
\[
 h(g_1,\ldots,g_m)=\alpha_1(g_1)\cdots\alpha_m(g_m).
\]
Theorem~\ref{thm:groups} applies to these endomorphisms.
\end{proof}
For a direct-product homomorphism, images of different $\alpha_j$
commute. Theorem~\ref{thm:groups} allows endomorphisms without that
property, including an ordered product of the unmodified inputs.

\begin{corollary}[Two-dimensional recurrences]\label{cor:plane}
Let $P\subset\Z^2$ be a finite nonempty set of lexicographically
positive predecessor vectors, and define $I$ and $C$ as above.
Let $G$ be a finite $p$-group, let $\alpha_1,\ldots,\alpha_m$
be endomorphisms of $G$, and suppose the boundary faces and the additional sequence are $p$-automatic, with compatible boundary values.
Then the unique array defined by the ordered-product recurrence in
Theorem~\ref{thm:groups} is $p$-automatic. The conclusion includes
homomorphisms from $G^m$ to $G$ and the additive abelian rules of
Theorem~\ref{thm:abelian}, without a further sign restriction on $P$.
\end{corollary}
The proof in Section~\ref{sec:planar-proof} checks the geometric and
boundary hypotheses and explains why all planar predecessor systems
in this statement are covered.

The number of coordinates containing a negative entry is not itself
the boundary difficulty. A coordinate with only nonpositive predecessor
entries has width $d_i=0$, hence no prescribed boundary face.
The following theorem allows such coordinates without restriction.

\begin{theorem}[At most one coordinate with both signs]\label{thm:mixed}
Let $P\subset\Z^n$ be finite, nonempty and lexicographically positive.
Suppose at most one coordinate $i$ satisfies
\[
 \min_j(v_j)_i<0<\max_j(v_j)_i.
\]
For any finite $p$-group $G$, ordered endomorphisms
$\alpha_1,\ldots,\alpha_m$, compatible automatic boundary faces,
and an automatic additional sequence, the unique solution of the recurrence in
Theorem~\ref{thm:groups} is $p$-automatic.
\end{theorem}
Section~\ref{sec:mixed-geometry} proves the theorem by removing
negative entries in coordinates without prescribed boundary faces.
The fixed obstruction in Corollary~\ref{cor:uniform-obstruction} has
two coordinates with both signs and fails for every prime. Thus the
threshold in the theorem cannot be increased using their number alone.

\section{Finite kernels and elementary constructions}\label{sec:elementary}

Fixing the low digits of the indices leaves another array to compute.
The $p$-kernel collects these remaining arrays. A finite kernel will give
the states of an automaton, and products of automata will let us combine
the prescribed data and test the boundary conditions.
The kernel characterization is treated in
\cite[Section~6.6]{AS}, with the multidimensional formulation in
\cite[Section~14.2]{AS}.

\subsection{Digit words and finite kernels}
\label{sec:digit-automata}

Write $\mathcal D_n=\{0,\ldots,p-1\}^n$ for the alphabet of digit
vectors. For $u:\N^n\to A$, where $A$ is finite, define
\[
 u_{e,\rct}(\x)=u(p^e\x+\rct),\qquad
 \kerp(u)=\{u_{e,\rct}:e\ge0,\ 0\le r_i<p^e\}.
\]
The array $u_{e,\rct}$ is called a section of $u$. Its indices have
their $e$ low digits fixed by $\rct$. For $e=0$, the only choice is
$\rct=0$, so $u_{0,0}=u$.

We write an automaton with output as
$\mathcal A=(Q,q_0,\delta,\tau)$. Here $Q$ is its finite state set,
$q_0$ is its initial state, $\delta:Q\times\mathcal D_n\to Q$ is
its transition map, and $\tau:Q\to A$ gives the output of a state.
For a word $w$ of digit vectors, the extended transition $\delta^*$
means that the letters of $w$ are read in their written order. Formally,
\[
 \delta^*(q,\varepsilon)=q,\qquad
 \delta^*(q,w\rct)=\delta(\delta^*(q,w),\rct),
\]
where $\varepsilon$ denotes the empty word.

The output on a word $w$ is
$\mathcal A(w)=\tau(\delta^*(q_0,w))$.
For $w=\rct_0\cdots\rct_{\ell-1}$, the represented index is
\[
 \begin{aligned}
 \sum_{j=0}^{\ell-1}p^{\ell-1-j}\rct_j
 &\quad\text{when most significant digits are read first},\\
 \sum_{j=0}^{\ell-1}p^j\rct_j
 &\quad\text{when least significant digits are read first}.
 \end{aligned}
\]
Each coordinate has the same word length. Adding a zero digit vector
at the most significant end leaves the index unchanged. Thus padding
invariance means $\mathcal A(\boldsymbol0w)=\mathcal A(w)$ in the
first convention and $\mathcal A(w\boldsymbol0)=\mathcal A(w)$ in
the second. The convention will be specified for each construction.

The empty word represents the zero tuple and has output $u(0)$.
For a machine already correct on all nonempty padded words, this
empty-word convention can be obtained by adding a new initial state
with output $u(0)$ and the same outgoing transitions as the old
initial state. Every nonempty word then has its previous output.
When the output alphabet is $\{0,1\}$, we also call the machine a
\emph{recognizer}. It accepts a word if its final output is $1$.
\par

\begin{lemma}[Finite-kernel characterization]\label{lem:kernel}
A finite-valued array $u$ is $p$-automatic if and only if $\kerp(u)$
is finite. When the kernel is finite, an automaton for $u$ has
\[
 Q=\kerp(u),\qquad q_0=u,\qquad
 \delta(v,\sct)(\x)=v(p\x+\sct),\qquad \tau(v)=v(0).
\]
It reads digits from least to most significant. The two digit-reading
conventions define the same class of arrays.
\end{lemma}
\begin{proof}
\textit{From a finite kernel to an automaton.}
First check that a transition remains in $Q$. If $v=u_{e,\rct}$, then
\[
 \delta(v,\sct)(\x)
 =u\bigl(p^{e+1}\x+\rct+p^e\sct\bigr)
 =u_{e+1,\rct+p^e\sct}(\x).
\]
Since $0\le r_i+p^es_i<p^{e+1}$, this is a member of $\kerp(u)$.

Let $\rct_0,\ldots,\rct_{\ell-1}$ be the digit vectors read so far,
starting with the least significant ones, and put
$\boldsymbol b_\ell=\sum_{j<\ell}p^j\rct_j$.
Induction on $\ell$ gives the state identity
\[
 \delta^*(u,\rct_0\cdots\rct_{\ell-1})(\x)
       =u(p^\ell\x+\boldsymbol b_\ell).
\]
The initial case is $u(\x)$, and the transition just computed proves
the induction step. At the end of the input, setting $\x=0$ gives
output $u(\boldsymbol b_\ell)$, the requested entry. Extra zero digit
vectors at the end do not change the output, because
$\delta(v,0)(0)=v(0)$.

\textit{From an automaton to a finite kernel.}
Suppose first that $\mathcal A=(Q,q_0,\delta,\tau)$ reads most
significant digits first. Let $w(\x)$ be a padded digit word for
$\x$, and put $q(\x)=\delta^*(q_0,w(\x))$.
For $0\le r_i<p^e$, let $\rho_e(\rct)$ be the word of exactly
$e$ digit vectors representing $\rct$ in that same order, including
initial zero vectors where necessary. Appending that word to $w(\x)$
represents $p^e\x+\rct$. Define an output map on the existing states by
\[
 \tau_{e,\rct}(q)=\tau\bigl(\delta^*(q,\rho_e(\rct))\bigr).
\]
Then
\[
 u_{e,\rct}(\x)=\tau_{e,\rct}(q(\x)).
\]
The map $\x\mapsto q(\x)$ is fixed. Only the function
$\tau_{e,\rct}:Q\to A$ varies, and there are at most $|A|^{|Q|}$
such functions. Thus there are at most that many kernel arrays.

\textit{Changing the digit-reading direction.}
For this equivalence, see also \cite[Sections~4.3 and~5.2]{AS}.
To reverse the reading direction, for a digit vector $\rct$ write
$\delta_{\rct}(q)=\delta(q,\rct)$.
Use functions $\eta:Q\to A$ as new states, start with $\eta=\tau$,
and define
\[
 \widehat\delta(\eta,\rct)=\eta\circ\delta_{\rct},\qquad
 \widehat\tau(\eta)=\eta(q_0).
\]
After the word $\sct_1\cdots\sct_\ell$, the new state is
$\tau\circ\delta_{\sct_1}\circ\cdots\circ\delta_{\sct_\ell}$.
Its output is the old automaton's output on the reversed word
$\sct_\ell\cdots\sct_1$. There are at most $|A|^{|Q|}$ new states,
and they can be enumerated by closing the initial map under the
displayed transitions. Reversal exchanges initial zero padding with
final zero padding, so the old padding invariance gives the required
invariance in the new direction.
\end{proof}

\subsection{Products and output codings}\label{sec:product-automata}

To compute several arrays at the same index, we run their automata
on the same input word. The construction below specifies how the
component states and outputs are combined. It is the usual product
construction, together with an output map
\cite[Sections~4.3 and~5.4]{AS}.

\begin{definition*}[Product automaton]
Let $t\ge1$, and let
$\mathcal A_j=(Q_j,q_j,\delta_j,\tau_j)$, for $1\le j\le t$,
be automata over the same input alphabet $\mathcal D_n$, using the
same digit-reading convention. Write $A_j$ for their output alphabets.
Their \emph{product automaton} is
$\mathcal A=(Q,q_0,\delta,\tau)$, with
\[
 Q=\prod_{j=1}^tQ_j,\qquad q_0=(q_1,\ldots,q_t),\qquad
 \tau:Q\longrightarrow\prod_{j=1}^t A_j,
\]
and, for $\rct\in\mathcal D_n$,
\begin{align*}
 \delta((s_1,\ldots,s_t),\rct)
   &=(\delta_1(s_1,\rct),\ldots,\delta_t(s_t,\rct)),\\
 \tau(s_1,\ldots,s_t)&=(\tau_1(s_1),\ldots,\tau_t(s_t)).
\end{align*}
\end{definition*}

A product state records one state from each component machine.
Induction on the length of $w$ gives
\begin{equation}\label{eq:product-state}
 \delta^*(q_0,w)=
 \bigl(\delta_1^*(q_1,w),\ldots,\delta_t^*(q_t,w)\bigr).
\end{equation}
Consequently, if $\mathcal A_j$ computes $u_j$, the product computes
$\x\mapsto(u_1(\x),\ldots,u_t(\x))$. Its state set has
$\prod_j|Q_j|$ elements. Padding preserves its output because it
preserves every component output. The states themselves may change
under padding. If the machines initially use different reading
directions, Lemma~\ref{lem:kernel} first converts them to a common one.

A \emph{coding} is a map between finite output alphabets, applied at
each index. For any finite set $B$ and map
$h:\prod_j A_j\to B$, replacing the product's output map $\tau$
by $h\circ\tau$ computes
\[
 \x\longmapsto h(u_1(\x),\ldots,u_t(\x)).
\]
The states, initial state and transitions are unchanged. For example,
if the $A_j$ are a common group $G$, the coding
$h(a_1,\ldots,a_t)=a_1\cdots a_t$ computes their ordered product
at each index. This applies to a nonabelian group because the order
of the factors is fixed. Boolean codings give equality tests and
unions, intersections and complements of recognized sets.

\paragraph{A two-dimensional example.}
In base $2$, consider $u(m,n)=m\bmod2$ and $v(m,n)=n\bmod3$.
Use most significant digits first, with input $(r,s)\in\{0,1\}^2$.
The first machine has states $\{0,1\}$, initial state $0$, output
equal to its state, and transition $a\mapsto r$. The second has
states $\{0,1,2\}$, initial state $0$, output equal to its state,
and transition $b\mapsto 2b+s\pmod3$. These transitions follow
from appending a binary digit to the corresponding coordinate.
The product therefore has six states and transition
\[
 (a,b)\xrightarrow{(r,s)}(r,2b+s\bmod3).
\]
For $(m,n)=(5,2)$, the padded binary representations are $101$ and
$010$, so the path is
\[
 (0,0)\xrightarrow{(1,0)}(1,0)
 \xrightarrow{(0,1)}(0,1)
 \xrightarrow{(1,0)}(1,2).
\]
The output $(1,2)$ is $(5\bmod2,2\bmod3)$.
Changing the output to $h(a,b)=1$ if $(a,b)=(0,0)$ and $0$ otherwise
makes the same product a recognizer for the indices where $m$ is
even and $n$ is divisible by $3$.

\paragraph{Reachable states and comparison of outputs.}
A state is \emph{reachable} if it equals $\delta^*(q_0,w)$ for some
finite input word $w$, including the empty word. Starting with
$R_0=\{q_0\}$, form
\[
 R_{h+1}=R_h\cup
 \{\delta(q,\rct):q\in R_h,\ \rct\in\mathcal D_n\}.
\]
This increasing sequence of subsets of the finite set $Q$ stabilizes.
Induction shows that $R_h$ contains exactly the states reached by
words of length at most $h$, so the stabilized set contains all
reachable states. It is closed under transitions. Restricting the
machine to this set preserves every output and may reduce the state
count below $\prod_j|Q_j|$.

To test whether two machines with the same output alphabet compute
the same array, give their product the Boolean output
$h(a,b)=\one_{\{a=b\}}$. The arrays agree at all indices if and
only if every reachable product state outputs $1$. A reachable state
with output $0$ supplies a counterexample: retain a word reaching it
and decode that word by the chosen digit convention. The same
procedure checks any specified relation among finitely many array
values. Section~\ref{sec:certified-construction} applies it to a
recurrence and its boundary conditions.

\paragraph{Merging equivalent states.}
For a fixed reading convention, two states are equivalent if they
give the same output after every continuation word. This can also
be checked by a finite computation. On the reachable state set, put
\[
 q\sim_0 s\quad\Longleftrightarrow\quad\tau(q)=\tau(s),
\]
and define $q\sim_{h+1}s$ to mean that $q\sim_hs$ and
$\delta(q,\rct)\sim_h\delta(s,\rct)$ for every input $\rct$.
Induction on $h$ shows that $q\sim_hs$ if and only if the outputs
agree after all words of length at most $h$. The partitions into
equivalence classes successively refine a finite state set, so they
stabilize. At stabilization, equivalent states have equal outputs
and equivalent successors. They can therefore be replaced by one
state per class, with the induced transitions and output.

This operation is called \emph{minimization}. Distinct classes are
separated by some continuation word. Choose a word reaching each
class from the initial state. Any machine with the same output on
all input words must reach different states on these chosen words,
since otherwise it could not distinguish their separating
continuations. Thus the quotient has the least possible number of
states in the chosen reading convention. This is the minimization
used for the counts in Section~\ref{sec:binary-certificates}.
\par

\subsection{Face tests, translations and boundary data}

\begin{lemma}[Closure operations]\label{lem:closure}
For any positive integer $t$, the tuple of $t$ $p$-automatic
arrays and any coding of that tuple are $p$-automatic.
Fixed-coordinate equality tests, finite unions
of boundary faces, and their complements are automatic. If $u$ takes
values in a finite set with distinguished value $0$, then, for fixed
$\vct\in\Z^n$, its translated zero extension
\[
 (T_{\vct}u)(\x)=
 \begin{cases}u(\x-\vct),&\x-\vct\in\N^n,\\
 0,&\text{otherwise}
 \end{cases}
\]
is automatic. All these constructions are effective.
\end{lemma}
\begin{proof}
\textit{Products and output codings.}
The product construction and the coding argument in
Section~\ref{sec:product-automata} prove the first assertion.
Their transition tables are obtained directly from those of the
input machines.

\textit{Face tests.}
Fix a coordinate $i$ and an integer $c\ge0$. To recognize $x_i=c$
from least significant digits first, start in state $c$ and use the
state set $\{0,\ldots,c\}\cup\{\bot\}$, where $\bot$ records a
mismatch. On a digit vector $\rct$, a state $a\ne\bot$ moves to
$\lfloor a/p\rfloor$ if $r_i$ is the remainder of $a$ modulo $p$,
and otherwise moves to $\bot$. The state $\bot$ is absorbing.
Output $1$ in state $0$ and $0$ in all other states.
After $\ell$ matching digits, the state is $\lfloor c/p^\ell\rfloor$,
the part of $c$ still to be checked. Thus the output is $1$ precisely
when the whole input equals $c$ in coordinate $i$.
The other coordinates are ignored. Products of these machines and
Boolean output functions give finite unions, intersections and
complements of the face tests.

\textit{Translations, including negative coordinates.}
Take a least-significant-digit automaton
$\mathcal A=(Q,q_0,\delta,\tau)$ for $u$, retaining its reachable
states. For $q\in Q$, let $u_q(\x)$ be the output obtained by starting
in state $q$ and reading $\x$. Its zero extension to $\Z^n$ is denoted
by $\widetilde u_q$. These arrays inherit padding invariance because
$q$ is reachable from the initial state. We will use states
$(q,\boldsymbol z)$ to represent
\[
 g_{q,\boldsymbol z}(\x)=\widetilde u_q(\x+\boldsymbol z).
\]
The initial state is $(q_0,-\vct)$, which represents $T_{\vct}u$.
For an input digit $\rct\in\mathcal D_n$, divide componentwise
\[
 \rct+\boldsymbol z=p\boldsymbol z'+\sct,
       \qquad 0\le s_i<p.
\]
The transition is
$(q,\boldsymbol z)\mapsto(\delta(q,\sct),\boldsymbol z')$.
Indeed,
\[
 g_{q,\boldsymbol z}(p\x+\rct)
 =\widetilde u_q\bigl(p(\x+\boldsymbol z')+\sct\bigr)
 =\widetilde u_{\delta(q,\sct)}(\x+\boldsymbol z').
\]
If a coordinate of $\x+\boldsymbol z'$ is negative, the same
coordinate of $p(\x+\boldsymbol z')+\sct$ is negative, because
$s_i<p$. Both sides are then zero, so the identity covers that case too.

The integer vectors $\boldsymbol z$ range over a finite set. Put
$L_i=\min(-v_i,0)$ and $U_i=\max(-v_i,0)$.
The initial vector belongs to $\prod_i[L_i,U_i]$, and
\[
 L_i\le\left\lfloor\frac{z_i+r_i}{p}\right\rfloor\le U_i
       \quad\text{when }L_i\le z_i\le U_i.
\]
For the lower bound, use $z_i+r_i\ge L_i$ and $L_i\le0$.
For the upper bound, use $z_i+r_i\le U_i+p-1$ and the fact that
$U_i$ is a nonnegative integer. Thus the displayed product of integer
intervals is invariant under transitions.
The output of $(q,\boldsymbol z)$ is
$\widetilde u_q(\boldsymbol z)$, evaluated by the given automaton
when $\boldsymbol z\in\N^n$, and set to zero otherwise.
The state set has at most $|Q|\prod_i(U_i-L_i+1)$ elements.
The state identities prove that the machine computes the translated
array. Its transition table and outputs are obtained by integer
arithmetic and evaluations of the given automaton at those states.
\end{proof}

\begin{lemma}[Boundary assembly]\label{lem:boundary}
If $c:C\to H$ is automatic on its faces and the face values agree,
its extension $c_0$ by the neutral element to $\N^n\setminus C$
is $p$-automatic and can be constructed effectively.
\end{lemma}
\begin{proof}
Write $e_H$ for the neutral element of $H$, using $0$ in additive
notation. Enumerate the faces as
\[
 F_\ell=\{\x:x_{i_\ell}=j_\ell\},\quad 1\le\ell\le s,
 \qquad s=\sum_{i=1}^n d_i.
\]
Since $P$ is nonempty and its vectors are lexicographically positive,
some $d_i$ is positive, so $s\ge1$.
Let $\pi_i:\N^n\to\N^{n-1}$ delete coordinate $i$, and let
$c_\ell:\N^{n-1}\to H$ be the prescribed array on $F_\ell$.
Define arrays on the whole orthant by
\[
 b_\ell(\x)=c_\ell(\pi_{i_\ell}\x),\qquad
 \chi_\ell(\x)=\one_{F_\ell}(\x).
\]
An automaton for $b_\ell$ ignores coordinate $i_\ell$ in each digit
vector and runs the automaton for $c_\ell$ on the remaining coordinates.
In dimension one, it simply outputs the prescribed value at the face's
single point. The indicator $\chi_\ell$ is automatic by
Lemma~\ref{lem:closure}.

To select a value when several faces meet, assign each point to the
first face containing it. The corresponding Boolean masks are
\[
 \mu_\ell(\x)=\chi_\ell(\x)
                   \prod_{h<\ell}(1-\chi_h(\x)).
\]
The masks use ordinary products of Boolean values, so they are
automatic $\{0,1\}$-valued arrays. At a point of $C$, one mask equals $1$, namely
the mask for the first containing face. Outside $C$, all masks vanish.

Define a coding on the full output alphabet by
\[
 \Phi((\varepsilon_1,\ldots,\varepsilon_s),(g_1,\ldots,g_s))=
 \begin{cases}
 g_{\min\{\ell:\varepsilon_\ell=1\}},&\text{if some }\varepsilon_\ell=1,\\
 e_H,&\text{if all }\varepsilon_\ell=0.
 \end{cases}
\]
Here each $\varepsilon_\ell$ belongs to $\{0,1\}$ and each $g_\ell$
belongs to $H$. Choosing the least index also defines $\Phi$ on
tuples that never occur as simultaneous machine outputs.
For the computed masks,
\[
 c_0(\x)=\Phi\bigl((\mu_\ell(\x))_{\ell=1}^s,
                       (b_\ell(\x))_{\ell=1}^s\bigr).
\]
The product automaton for all the masks and face arrays, with its
output composed with $\Phi$, therefore computes $c_0$.
\par
If $\x$ lies on several faces, compatibility says that the values
$b_\ell(\x)$ on those faces coincide. Hence the selected value agrees
with every prescribed restriction.
\end{proof}

For example, suppose $n=2$ and $d_1=d_2=1$. Write $c_1(y)$ for
the data on $x=0$ and $c_2(x)$ for the data on $y=0$.
Compatibility means $c_1(0)=c_2(0)$. Choosing $x=0$ first gives
\[
 c_0(x,y)=
 \begin{cases}
 c_1(y),&x=0,\\
 c_2(x),&x>0,\ y=0,\\
 e_H,&x>0,\ y>0.
 \end{cases}
\]
The choice of face order has no effect at the corner because the
prescribed values there agree.

\subsection{Existence and uniqueness}

\begin{lemma}[Existence and uniqueness]\label{lem:unique}
Let $P=\{\vct_1,\ldots,\vct_m\}$ be a finite nonempty set of
lexicographically positive vectors, with interior $I$ and boundary $C$
as defined in Section~\ref{sec:introduction}. Let $A$ be a nonempty
set. Any boundary function $c:C\to A$ and maps $h_{\x}:A^m\to A$
for $\x\in I$ determine a unique array $a:\N^n\to A$ satisfying
\[
 a(\x)=c(\x)\quad(\x\in C),\qquad
 a(\x)=h_{\x}\bigl(a(\x-\vct_1),\ldots,a(\x-\vct_m)\bigr)
       \quad(\x\in I).
\]
\end{lemma}
\begin{proof}
We assign an integer weight to each point so that every predecessor
has smaller weight. For vectors $w,x\in\Z^n$, the scalar product is
$w\cdot x=\sum_{i=1}^n w_i x_i$. Choose the positive integer weights
backwards by
\[
 w_n=1,\qquad
 w_i=1+\max_{1\le j\le m}\sum_{k=i+1}^n w_k|(v_j)_k|
       \quad(i=n-1,\ldots,1),
\]
and put $w=(w_1,\ldots,w_n)$ and $\rho(x)=w\cdot x$.
The maximum is over the $m$ predecessor vectors. If the first nonzero
entry of $\vct_j$ is in coordinate $i$, that entry is at least $1$.
Therefore
\[
 \rho(\vct_j)
 =w_i(v_j)_i+\sum_{k>i}w_k(v_j)_k
 \ge w_i-\sum_{k>i}w_k|(v_j)_k|\ge1.
\]
For $i=n$, the last inequality follows from $w_n=1$.

For each integer $N\ge0$, define
\[
 S_N=\{\x\in\N^n:\rho(\x)\le N\},\qquad S_{-1}=\varnothing.
\]
The set $S_N$ contains finitely many points, because its coordinates
satisfy $0\le x_i\le\lfloor N/w_i\rfloor$.
Induct on $N$ to assign the values on $S_N$, starting with the empty
assignment on $S_{-1}$. Suppose the values on $S_{N-1}$ are already
assigned. For each $\x\in S_N\setminus S_{N-1}$, prescribe $c(\x)$
if $\x\in C$. Otherwise $\x\in I$, so every $\x-\vct_j$ belongs
to $\N^n$ and
\[
 \rho(\x-\vct_j)=N-\rho(\vct_j)\le N-1.
\]
All inputs to $h_{\x}$ have thus been assigned, and the rule gives
$a(\x)$. Points of weight $N$ do not depend on one another. If no
point has that weight, then $S_N=S_{N-1}$ and the step adds no values.
The induction does not require every integer to occur as a weight.

Every point of $\N^n$ belongs to some $S_N$, so these assignments
define a solution on the whole domain. If two solutions agree on
$S_{N-1}$, their prescribed boundary values and their identical inputs
to each $h_{\x}$ make them agree on $S_N$. Induction proves uniqueness.
\end{proof}

\subsection{Adjugates and division by a determinant}\label{sec:adjugate}

The generating functions below satisfy matrix equations over commutative
rings. To solve those equations, we need a determinant identity that
remains valid when the coefficient ring has zero divisors.

\begin{definition}
For an $s\times s$ matrix $D$ over a commutative ring, let
$D^{(i,j)}$ be the matrix obtained by deleting row $i$ and column $j$.
The \emph{adjugate} of $D$ is the transpose of its cofactor matrix,
so
\[
 (\adj D)_{ij}=(-1)^{i+j}\det D^{(j,i)}.
\]
The determinant of the empty matrix is $1$. In particular, the adjugate
of a $1\times1$ matrix is the $1\times1$ identity matrix.
\end{definition}

\begin{lemma}[Adjugate identity]\label{lem:adjugate}
For a square matrix $D$ over a commutative ring $A$,
\[
 D\adj D=(\det D)I_s=(\adj D)D.
\]
If $\det D$ is invertible in $A$, then
$D^{-1}=(\det D)^{-1}\adj D$. If $A$ is an integral domain and
$\det D\ne0$, the same inverse formula holds over its fraction field.
\end{lemma}
\begin{proof}
The $(i,j)$ entry of $D\adj D$ is
$\sum_k d_{ik}(-1)^{j+k}\det D^{(j,k)}$.
For $i=j$, this is the expansion of $\det D$ along row $j$.
For $i\ne j$, it is the determinant of the matrix obtained by replacing
row $j$ by row $i$, so it vanishes because two rows coincide.
Expansion by columns proves the other identity. These determinant
identities use only sums and products in a commutative ring.
Multiplication by $(\det D)^{-1}$ proves the inverse statements.
\end{proof}

For variables $\X=(X_1,\ldots,X_n)$, the notation $A[[\X]]$
denotes the ring of formal power series with coefficients in $A$.
Its sums and products are defined coefficient by coefficient, without
a convergence requirement.

\begin{lemma}[Units in a power-series ring]\label{lem:seriesunit}
A series $Q\in A[[\X]]$ over a commutative ring $A$ is invertible
if and only if $Q(0)$ is invertible in $A$.
Consequently, if a polynomial matrix $D$ satisfies $D(0)=I_s$,
then $\det D$ and $D$ are invertible over $A[[\X]]$.
\end{lemma}
\begin{proof}
The constant coefficient of $QV=1$ gives $Q(0)V(0)=1$, proving
necessity. Conversely, if $c=Q(0)$ is invertible, write $Q=c(1-W)$
with $W(0)=0$. The series $c^{-1}\sum_{r\ge0}W^r$ is well-defined
because a coefficient of total degree $d$ receives terms only from
$r\le d$. Multiplication by $Q$ gives $1$.
If $D(0)=I_s$, then $(\det D)(0)=1$, and the adjugate identity
gives the matrix inverse.
\end{proof}

In Section~\ref{sec:nonnegative}, the determinant has constant
coefficient $1$, so the inverse is a matrix of formal series even over
$\Z/p^e\Z$. In Section~\ref{sec:negative}, a determinant with zero
constant coefficient is divided out in a fraction field. The argument
there separately establishes that the resulting vector is a formal
power series.

\subsection{The one-dimensional case}\label{sec:one-dimension}

When $n=1$, a lexicographically positive predecessor vector is a
positive integer. Write the predecessor offsets as $r_1,\ldots,r_m$
and put $d=\max_j r_j$. Then the prescribed boundary is the initial
segment $C=\{0,\ldots,d-1\}$. Each boundary face is a single point,
so any assignment of these $d$ values meets the boundary hypothesis.

\begin{corollary}[One-dimensional recurrences]\label{cor:one-dimension}
Let $G$ be a finite $p$-group, let $r_1,\ldots,r_m$ be positive
integers, and let $\alpha_1,\ldots,\alpha_m$ be endomorphisms of $G$.
For any prescribed $a(0),\ldots,a(d-1)\in G$ and any $p$-automatic
additional sequence $b:\N\to G$, the recurrence
\[
 a(k)=\alpha_1(a(k-r_1))\cdots\alpha_m(a(k-r_m))b(k)
       \qquad(k\ge d)
\]
has a unique solution, and this solution is $p$-automatic.
\end{corollary}
\begin{proof}
Lemma~\ref{lem:unique} applies with the single weight $w_1=1$.
All offsets are positive, so the geometric hypothesis of
Theorem~\ref{thm:groups} holds. The boundary has $d$ points and is
automatic on its faces. That theorem therefore gives automaticity.
Its proof uses the nonnegative construction in this case, with no
algebraic division argument for negative predecessor entries.
Section~\ref{sec:one-dimensional-division} writes out a scalar instance
of the construction.
\end{proof}

If $b(k)=\beta$ is constant, eventual periodicity follows directly
from the finite number of consecutive blocks. For $k\ge d-1$, put
\[
 s_k=(a(k-d+1),\ldots,a(k))\in G^d.
\]
For $s=(g_0,\ldots,g_{d-1})$, the next entry is
\[
 h(s)=\alpha_1(g_{d-r_1})\cdots\alpha_m(g_{d-r_m})\beta,
\]
and the next block is $(g_1,\ldots,g_{d-1},h(s))$. For $d=1$,
this block consists just of $h(s)$. Thus a fixed map $G^d\to G^d$
updates the blocks. There are $|G|^d$ possible blocks, so two blocks
in the first $|G|^d+1$ positions coincide. All subsequent blocks
repeat with the same period because the update map is fixed. Reading
their last coordinates proves eventual periodicity of $a$.
An eventually periodic sequence is automatic in every integer base
at least $2$. An automaton can test its finite initial segment and
compute the index modulo its period for the remaining entries.

For a general automatic $b$, the next block also depends on $b(k+1)$.
Repeated blocks alone then give no periodicity conclusion. The
automaticity assertion in Corollary~\ref{cor:one-dimension} uses the
group rule and the automaticity of $b$.

\section{Automatic series divided by a polynomial}\label{sec:division}

An $n$-dimensional array and its generating series in $n$ variables
contain the same coefficients. Section operators select coefficients
according to the low digits of their indices. In one variable over
$\F_p$, these are the operators used by Christol, Kamae,
Mend\`es France and Rauzy \cite[Section~6, p.~407]{CKMR}, who work
more generally over $\F_q$ in base $q$.
We use the name \emph{Cartier operators}, as in
\cite[Definition~1.7]{RY15}, and explain its relation to Cartier's
operator on differential forms after the definition. The product
identity will let us compute sections while keeping a fixed denominator.

Fix $R=\Z/p^e\Z$ with $e\ge1$, and retain the digit alphabet
$\mathcal D_n=\{0,\ldots,p-1\}^n$. For
$\X=(X_1,\ldots,X_n)$ and $\boldsymbol a\in\N^n$, write
$\X^{\boldsymbol a}=\prod_iX_i^{a_i}$.
A formal series has the form
\[
 F(\X)=\sum_{\boldsymbol a\in\N^n}f_{\boldsymbol a}\X^{\boldsymbol a},
       \qquad f_{\boldsymbol a}\in R.
\]
We write $[\X^{\boldsymbol a}]F=f_{\boldsymbol a}$ for coefficient
extraction. No analytic convergence is involved. In a product of formal
series, each coefficient is a finite sum.

\begin{definition}[Cartier operator]\label{def:cartier}
For $\rct\in\mathcal D_n$, the Cartier operator $\Lambda_{\rct}$
is the $R$-linear map given by
\[
 \Lambda_{\rct}F
   =\sum_{\x\in\N^n}f_{p\x+\rct}\X^{\x}.
\]
It selects the coefficients whose index vector has low digit $\rct$
and divides the remaining indices by $p$ after removing that digit.
Its coefficient array is the section $\x\mapsto f_{p\x+\rct}$.
On a monomial the rule is
\[
 \Lambda_{\rct}(\X^{\boldsymbol a})=
 \begin{cases}
 \X^{(\boldsymbol a-\rct)/p},&\boldsymbol a\equiv\rct\pmod p,\\
 0,&\text{otherwise},
 \end{cases}
\]
with the congruence interpreted coordinatewise.
\end{definition}

\paragraph{Relation to Cartier's operator.}\label{sec:cartier-attribution}
Cartier \cite[p.~427, formulas~(5)]{Cartier} introduces an operator on
differential forms in characteristic $p$. In the field $\F_p((X))$
of formal Laurent series, write a differential form as
$\omega=F(X)\,dX$ and set $dH=H'(X)\,dX$ for $F,H\in\F_p((X))$.
His additive operator $\mathcal C$ satisfies
\[
 \mathcal C(H^p\omega)=H\mathcal C(\omega),\qquad
 \mathcal C(dH)=0,\qquad
 \mathcal C(X^{p-1}\,dX)=dX.
\]
For $0\le r<p-1$, the form $X^r\,dX$ is the differential of
$X^{r+1}/(r+1)$, so its image is zero. Decomposing a series over
$\F_p$ as $F=\sum_{r=0}^{p-1}X^r(\Lambda_rF)^p$ gives the first
identity below. Applying it to $X^{p-1-r}F$ gives the second for
$0\le r<p$:
\[
 \mathcal C(F\,dX)=(\Lambda_{p-1}F)\,dX,\qquad
 \mathcal C(X^{p-1-r}F\,dX)=(\Lambda_rF)\,dX.
\]
These identities explain the association with Cartier. Our definition
selects coefficients directly, which also makes sense in several
variables and over $\Z/p^e\Z$. The arguments below use this
coefficient definition and the identities proved next.

For example, in one variable and base $2$, the two operators are
\begin{align*}
 \Lambda_0(f_0+f_1X+f_2X^2+\cdots)&=f_0+f_2X+f_4X^2+\cdots,\\
 \Lambda_1(f_0+f_1X+f_2X^2+\cdots)&=f_1+f_3X+f_5X^2+\cdots.
\end{align*}
Thus an even-index coefficient and an odd-index coefficient move into
different series. The next identity reconstructs $F$ from those series.
The second identity explains how to select coefficients of a product
when one factor uses only powers divisible by $p$.

\begin{lemma}[Coefficient selection identities]\label{lem:cartier}
The Cartier operators are $R$-linear. For formal series $F,V,W$,
\begin{align}
 F(\X)&=\sum_{\rct\in\mathcal D_n}\X^{\rct}
                (\Lambda_{\rct}F)(\X^p),\label{eq:sectionsplit}\\
 \Lambda_{\rct}\bigl(V(\X)W(\X^p)\bigr)
       &=(\Lambda_{\rct}V)(\X)W(\X).\label{eq:sectionprod}
\end{align}
Here $\X^p=(X_1^p,\ldots,X_n^p)$.
If $V$ is a polynomial with $\deg_{X_i}V\le d_i$, then
$\deg_{X_i}\Lambda_{\rct}V\le\lfloor d_i/p\rfloor$
whenever the selected polynomial is nonzero.
\end{lemma}
\begin{proof}
For $a,b\in R$ and formal series $V,W$, coefficient selection gives
\[
 [\X^{\x}]\Lambda_{\rct}(aV+bW)
 =a[\X^{p\x+\rct}]V+b[\X^{p\x+\rct}]W.
\]
This proves linearity. For reconstruction, divide each coordinate of
an exponent vector by $p$. There is a unique expression
$\boldsymbol a=p\x+\rct$ with $\x\in\N^n$ and
$\rct\in\mathcal D_n$. Grouping the terms of $F$ by this remainder gives
\[
 F(\X)
 =\sum_{\rct\in\mathcal D_n}\sum_{\x\in\N^n}
      f_{p\x+\rct}\X^{p\x+\rct}
 =\sum_{\rct\in\mathcal D_n}\X^{\rct}
      (\Lambda_{\rct}F)(\X^p).
\]
Each monomial appears in one summand, which justifies the identity
coefficient by coefficient.

For the product identity, write
$V=\sum_{\boldsymbol a}v_{\boldsymbol a}\X^{\boldsymbol a}$ and
$W=\sum_{\boldsymbol b}w_{\boldsymbol b}\X^{\boldsymbol b}$.
For a fixed $\x\in\N^n$,
\begin{align*}
 [\X^{\x}]\Lambda_{\rct}\bigl(V(\X)W(\X^p)\bigr)
 &= [\X^{p\x+\rct}]\bigl(V(\X)W(\X^p)\bigr)\\
 &=\sum_{\boldsymbol a+p\boldsymbol b=p\x+\rct}
            v_{\boldsymbol a}w_{\boldsymbol b}\\
 &=\sum_{0\le\boldsymbol b\le\x}
            v_{p(\x-\boldsymbol b)+\rct}w_{\boldsymbol b}.
\end{align*}
In the middle sum both exponent vectors are nonnegative. Since
$0\le r_i<p$, its index equation forces $b_i\le x_i$ and then
$a_i=p(x_i-b_i)+r_i$. The last sum is the coefficient of
$\X^{\x}$ in $(\Lambda_{\rct}V)W$, which proves
\eqref{eq:sectionprod}.

Finally, if a monomial of $V$ survives selection, its exponent in
coordinate $i$ has the form $px_i+r_i\le d_i$.
Its new exponent therefore satisfies
$x_i\le\lfloor(d_i-r_i)/p\rfloor\le\lfloor d_i/p\rfloor$.
This proves the degree bound.
\end{proof}

To obtain an automaton, we need a finite family of series that remains
inside the family after each digit selection. The series themselves
will be the states. The output of a state will be its constant
coefficient, because reading the remaining digits eventually reduces
the requested index to zero.

\begin{proposition}[Cartier construction of an automaton]\label{prop:cartier}
Suppose a finite set $\mathcal S\subset R[[\X]]$ contains $F$ and
satisfies $\Lambda_{\rct}V\in\mathcal S$ for all
$V\in\mathcal S$ and $\rct\in\mathcal D_n$.
Define
\[
 \mathcal A_F=(\mathcal S,F,\delta,\tau),\qquad
 \delta(V,\rct)=\Lambda_{\rct}V,\qquad \tau(V)=[\X^0]V.
\]
Then $\mathcal A_F$ computes the coefficient array of $F$ by reading
coordinate digits together from least to most significant.
\end{proposition}
\begin{proof}
Write $F=\sum f_{\x}\X^{\x}$. After reading the first $\ell$
digit vectors $\rct_0,\ldots,\rct_{\ell-1}$, put
\[
 V_0=F,\qquad V_{\ell+1}=\Lambda_{\rct_\ell}V_\ell,
 \qquad \boldsymbol b_\ell=\sum_{j<\ell}p^j\rct_j.
\]
The state remains in $\mathcal S$ by the closure assumption.
We prove the coefficient identity
\begin{equation}\label{eq:cartier-invariant}
 [\X^{\x}]V_\ell=f_{p^\ell\x+\boldsymbol b_\ell}
       \qquad(\x\in\N^n).
\end{equation}
For $\ell=0$, it is the definition of $F$.
If it holds for $\ell$, then
\begin{align*}
 [\X^{\x}]V_{\ell+1}
 &= [\X^{p\x+\rct_\ell}]V_\ell\\
 &= f_{p^\ell(p\x+\rct_\ell)+\boldsymbol b_\ell}\\
 &= f_{p^{\ell+1}\x+\boldsymbol b_{\ell+1}},
\end{align*}
which proves the induction step.
If the full input has length $L$ and represents $\boldsymbol m$, then
$\boldsymbol b_L=\boldsymbol m$. Equation~\eqref{eq:cartier-invariant}
at $\x=0$ gives
$\tau(V_L)=f_{\boldsymbol m}$.
Appending a zero digit vector preserves the output because
$[\X^0]\Lambda_0V=[\X^0]V$.
\end{proof}

This is the machine of Lemma~\ref{lem:kernel} expressed in terms of
generating series. To implement it, one needs a finite list of state
representatives, their constant coefficients, and an exact identity
for each transition $\Lambda_{\rct}V$. Those identities give the
complete transition table. Starting from $F$, one may discard states
that are never reached. The next example computes the list and all
transitions directly. Theorem~\ref{thm:division} represents
states by polynomial tuples and bounds the number of such tuples.

\subsection{A complete example: Pascal's triangle modulo 2}
\label{sec:pascal-cartier}

We construct the automaton for the array
\[
 b(n,k)=\binom nk\bmod2,\qquad n,k\in\N,
\]
where $\binom nk=0$ for $k>n$. The first coordinate $n$ is the row
number and the second coordinate $k$ is the column number.
The boundary values are $b(n,0)=1$ for all $n$ and $b(0,k)=0$
for $k>0$. Pascal's identity gives
\[
 b(n,k)=b(n-1,k)+b(n-1,k-1)\quad(n,k\ge1)
\]
in $\F_2$. Table~\ref{tab:pascal-values} displays the first eight
rows, including the zeros beyond the end of each row.

\begin{table}[!htbp]
\centering
\begin{tabular}{c|cccccccc}
$n\backslash k$ & $0$ & $1$ & $2$ & $3$ & $4$ & $5$ & $6$ & $7$\\
\hline
$0$ & $1$ & $0$ & $0$ & $0$ & $0$ & $0$ & $0$ & $0$\\
$1$ & $1$ & $1$ & $0$ & $0$ & $0$ & $0$ & $0$ & $0$\\
$2$ & $1$ & $0$ & $1$ & $0$ & $0$ & $0$ & $0$ & $0$\\
$3$ & $1$ & $1$ & $1$ & $1$ & $0$ & $0$ & $0$ & $0$\\
$4$ & $1$ & $0$ & $0$ & $0$ & $1$ & $0$ & $0$ & $0$\\
$5$ & $1$ & $1$ & $0$ & $0$ & $1$ & $1$ & $0$ & $0$\\
$6$ & $1$ & $0$ & $1$ & $0$ & $1$ & $0$ & $1$ & $0$\\
$7$ & $1$ & $1$ & $1$ & $1$ & $1$ & $1$ & $1$ & $1$\\
\end{tabular}
\caption{The values $b(n,k)=\binom nk\bmod2$ for $0\le n,k<8$.
Rows are indexed by $n$ and columns by $k$. Entries with $k>n$ are zero.
The automaton in Table~\ref{tab:pascal-machine} computes the whole
infinite array.}
\label{tab:pascal-values}
\end{table}

\paragraph{The generating function.}
Work in $\F_2[[X,Y]]$. By the binomial theorem, the polynomial for
row $n$ is $\sum_{k\ge0}b(n,k)Y^k=(1+Y)^n$.
Thus the generating series is
\[
 F(X,Y)=\sum_{n,k\ge0}b(n,k)X^nY^k
       =\sum_{n\ge0}\bigl(X(1+Y)\bigr)^n
       =\frac1{1+X+XY}.
\]
The geometric sum is valid formally because $X(1+Y)$ has zero
constant coefficient. Write $Q(X,Y)=1+X+XY$, so $F=1/Q$.
The constant coefficient of $Q$ is $1$, which makes this division
valid by Lemma~\ref{lem:seriesunit}.

\paragraph{All four Cartier images.}
There are four input digits, $(r,s)\in\{0,1\}^2$, where $r$ is
a digit of $n$ and $s$ is a digit of $k$.
In characteristic $2$, squaring gives
\[
 Q(X,Y)^2=1+X^2+X^2Y^2=Q(X^2,Y^2).
\]
It follows that
\[
 F(X,Y)=\frac{Q(X,Y)}{Q(X^2,Y^2)}
       =Q(X,Y)F(X^2,Y^2).
\]
Apply \eqref{eq:sectionprod} with $V=Q$ and $W=F$ to obtain
\[
 \Lambda_{r,s}F=(\Lambda_{r,s}Q)F.
\]
The monomials $1$, $X$ and $XY$ of $Q$ have remainder vectors
$(0,0)$, $(1,0)$ and $(1,1)$, respectively. Their exponents become
$(0,0)$ after selection. There is no monomial with remainder $(0,1)$.
Consequently,
\[
 \Lambda_{0,0}Q=1,\qquad \Lambda_{0,1}Q=0,\qquad
 \Lambda_{1,0}Q=1,\qquad \Lambda_{1,1}Q=1.
\]
Therefore the four images of $F$ are $F,0,F,F$, in that order.
All four images of the zero series are zero.
The family $\mathcal S=\{F,0\}$ is closed under the Cartier
operators, so Proposition~\ref{prop:cartier} applies.

\paragraph{The complete automaton.}
Label the state representing $F$ by $q_F$ and the state representing
the zero series by $q_Z$. The initial state is $q_F$, and the outputs
are $\tau(q_F)=F(0,0)=1$ and $\tau(q_Z)=0$.
Table~\ref{tab:pascal-machine} lists every transition. In particular,
once the digit pair $(0,1)$ is read in state $q_F$, the machine enters
$q_Z$ and remains there.

\begin{table}[!htbp]
\centering
\begin{tabular}{c|cccc|c}
State & $(0,0)$ & $(0,1)$ & $(1,0)$ & $(1,1)$ & Output\\
\hline
$q_F$ & $q_F$ & $q_Z$ & $q_F$ & $q_F$ & $1$\\
$q_Z$ & $q_Z$ & $q_Z$ & $q_Z$ & $q_Z$ & $0$\\
\end{tabular}
\caption{The two-state Cartier automaton for $\binom nk\bmod2$.
The initial state is $q_F$. The first and second components of each
input pair are binary digits of $n$ and $k$, respectively. Pairs are
read from least to most significant. The four middle columns give successor states.
The last column gives the output when the input ends.}
\label{tab:pascal-machine}
\end{table}

\begin{samepage}
For example, $5=(101)_2$ and $1=(001)_2$. Reading from the low digits
gives the pairs $(1,1),(0,0),(1,0)$ and the path
\[
 q_F\xrightarrow{(1,1)}q_F\xrightarrow{(0,0)}q_F
       \xrightarrow{(1,0)}q_F.
\]
The output is $1$, in agreement with $\binom51=5$ modulo $2$.
For $5=(101)_2$ and $2=(010)_2$, the path is
\[
 q_F\xrightarrow{(1,0)}q_F\xrightarrow{(0,1)}q_Z
       \xrightarrow{(1,0)}q_Z.
\]
Its output is $0$, in agreement with $\binom52=10$ modulo $2$.
Additional digit pairs $(0,0)$ preserve the output in both states.
\end{samepage}

\paragraph{What the states mean for the whole array.}
Equating coefficients in the four Cartier identities gives
\begin{align*}
 b(2n,2k)&=b(n,k),& b(2n,2k+1)&=0,\\
 b(2n+1,2k)&=b(n,k),& b(2n+1,2k+1)&=b(n,k).
\end{align*}
These identities hold for all $n,k\ge0$, including positions above
the triangle. If $n=\sum_j2^jn_j$ and $k=\sum_j2^jk_j$, the
transition table therefore says
\[
 b(n,k)=1\quad\Longleftrightarrow\quad k_j\le n_j
       \text{ for every }j.
\]
This follows from the proved identities for formal series and applies
to all indices. The two states are both reachable and have different
outputs, so a one-state automaton cannot compute this array.

\subsection{Finite polynomial states for division}\label{sec:polynomial-states}

The denominator argument below is the Cartier construction in
\cite[Proposition~1.9 and the proof of Theorem~2.1]{RY15}, in the form
of \cite[Theorem~5.1]{RY15} that keeps all variables, adapted to a
tuple of automatic numerator series. We will use the reconstruction identity to
separate their low digits and the product identity to preserve a common
denominator. Lemma~\ref{lem:seriesunit} justifies division by any
polynomial with constant coefficient $1$ over $R[[\X]]$.

\begin{lemma}\label{lem:frobenius}
For any $Q\in R[\X]$ and $E=p^{e-1}$,
\[
 Q(\X)^{pE}=Q(\X^p)^E\quad\text{in }R[\X].
\]
\end{lemma}
\begin{proof}
Lift coefficients to integers. Modulo $p$, the binomial theorem and
Fermat's congruence give $Q(\X)^p\equiv Q(\X^p)$. If $A\equiv B$
modulo $p^j$, $j\ge1$, the binomial expansion of $(B+p^jC)^p$ gives
$A^p\equiv B^p$ modulo $p^{j+1}$. Iterating this $e-1$ times gives
the displayed equality modulo $p^e$.
\end{proof}

\newpage
\begin{theorem}[Explicit division construction]\label{thm:division}
Let $U\in R[[\X]]$ have $M$ distinct arrays in its $p$-kernel, and let
$Q\in R[\X]$ have $Q(0)=1$. Put
\[
 E=p^{e-1},\qquad d_i=\deg_{X_i}Q,\qquad
 V=\prod_{i=1}^n(Ed_i+1).
\]
Then $U/Q$ is $p$-automatic and has a least-significant-digit automaton
with at most
\begin{equation}\label{eq:bound}
 |R|^{MV}=p^{eMV}
\end{equation}
states. The transitions and outputs can be constructed from the
kernel transitions of $U$ and the coefficients of $Q$.
\end{theorem}
\begin{proof}
\textit{State representatives and the initial state.}
List the kernel series as $U_1,\ldots,U_M$, with $U_1=U$.
Record their transitions by the index map $\sigma$, defined by
\[
 \Lambda_{\sct}U_j=U_{\sigma(j,\sct)}
       \quad(1\le j\le M,\ \sct\in\mathcal D_n).
\]
Let $\mathcal P$ be the finite set of polynomial tuples
$\boldsymbol P=(P_1,\ldots,P_M)$ such that
$\deg_{X_i}P_j\le Ed_i$ for all $i,j$.
A tuple represents the formal series
\begin{equation}\label{eq:representation}
 F_{\boldsymbol P}=\frac{\sum_{j=1}^M P_jU_j}{Q^E}.
\end{equation}
The initial tuple is
$\boldsymbol P^{\mathrm{init}}=(Q^{E-1},0,\ldots,0)$.
It satisfies the degree bounds and represents $U/Q$.

\textit{Computing a transition.}
Fix an input digit $\rct\in\mathcal D_n$. We seek a tuple
$\boldsymbol P'$ representing $\Lambda_{\rct}F_{\boldsymbol P}$.
Put $B=Q^{(p-1)E}$. Lemma~\ref{lem:frobenius} gives
$Q^EB=Q(\X^p)^E$, while \eqref{eq:sectionsplit} gives
\[
 U_j(\X)=\sum_{\sct\in\mathcal D_n}\X^{\sct}
                         U_{\sigma(j,\sct)}(\X^p).
\]
Multiply numerator and denominator of \eqref{eq:representation} by
$B$, and substitute these expressions for the $U_j$. Then
\begin{align*}
 F_{\boldsymbol P}
 &=\frac{\sum_{j,\sct}\X^{\sct}P_jB\,
                         U_{\sigma(j,\sct)}(\X^p)}{Q(\X^p)^E}\\
 &=\sum_{j,\sct}\X^{\sct}P_jB
             \left(\frac{U_{\sigma(j,\sct)}}{Q^E}\right)(\X^p).
\end{align*}
Here and below, $j$ ranges over $1,\ldots,M$ and $\sct$ over
$\mathcal D_n$. The formal quotients exist because $Q(0)=1$.
Applying \eqref{eq:sectionprod} to each summand yields
\[
 \Lambda_{\rct}F_{\boldsymbol P}
 =\frac{\sum_{j,\sct}\Lambda_{\rct}(\X^{\sct}P_jB)
                            U_{\sigma(j,\sct)}}{Q^E}.
\]
Group the terms multiplying the same series $U_\ell$. The resulting
transition is $\delta(\boldsymbol P,\rct)=\boldsymbol P'$, where
\begin{equation}\label{eq:transition}
 P'_\ell=\sum_{\substack{j,\sct\\\sigma(j,\sct)=\ell}}
                \Lambda_{\rct}(\X^{\sct}P_jB).
\end{equation}
Thus $F_{\boldsymbol P'}=\Lambda_{\rct}F_{\boldsymbol P}$.

\textit{Why the state set is finite and closed.}
The degree of $\X^{\sct}P_jB$ in $X_i$ is at most
\[
 (p-1)+Ed_i+(p-1)Ed_i=pEd_i+p-1.
\]
Lemma~\ref{lem:cartier} bounds the degree after selection by
$\lfloor(pEd_i+p-1)/p\rfloor=Ed_i$.
Summing the selected polynomials does not increase this bound, so
$\boldsymbol P'$ belongs to $\mathcal P$.
There are $V=\prod_i(Ed_i+1)$ possible monomials in each component
$P_j$, each with a coefficient in $R$. Therefore
$|\mathcal P|=|R|^{MV}$.

\textit{Outputs and termination of the construction.}
Assign the output
\[
 \tau(\boldsymbol P)=[\X^0]F_{\boldsymbol P}
             =\sum_{j=1}^M P_j(0)U_j(0).
\]
The last equality uses $Q(0)=1$.
Starting with $\boldsymbol P^{\mathrm{init}}$, compute its successors
for all $p^n$ input digits using \eqref{eq:transition}. Repeat for
new tuples, identifying tuples by their polynomial coefficient lists.
The process adds at most $|R|^{MV}$ tuples and therefore terminates.

Different tuples may represent the same series. The transition identity
still gives the coefficient invariant \eqref{eq:cartier-invariant}
for the represented series, so the proof of
Proposition~\ref{prop:cartier} gives the required outputs at all indices.
The bound counts the polynomial tuples before any further minimization.
\end{proof}

\begin{remark}[Obtaining the numerator's kernel]
Given a most-significant-digit automaton for $U$, enumerate the output
maps produced by all finite suffixes as in Lemma~\ref{lem:kernel}.
At most $|R|^{|Q|}$ distinct output maps can occur, where
$Q$ is the given state set. Each such map represents a section of
$U$, and composition with a digit transition gives its successor.
Redundant
states may be merged by ordinary finite-automaton minimization. The
construction uses the given transition table and output map
to enumerate all reachable maps. Bound
\eqref{eq:bound} also holds with the number of states in any such
section presentation in place of the exact kernel size $M$.
\end{remark}

\subsection{A one-dimensional affine recurrence}\label{sec:one-dimensional-division}

To apply the construction in one dimension, take $R=\Z/p^e\Z$,
a $p$-automatic sequence $b:\N\to R$, and the recurrence
\[
 a(0)=c,\qquad a(k)=a(k-1)+b(k)\quad(k\ge1).
\]
Let $A(X)=\sum_{k\ge0}a(k)X^k$ and
$B(X)=\sum_{k\ge0}b(k)X^k$. Summing the recurrence for $k\ge1$
gives $A-c=XA+B-b(0)$, hence
\[
 A(X)=\frac{c+B(X)-b(0)}{1-X}.
\]
The numerator has automatic coefficients. It is obtained from $b$
by replacing its value at $0$ with $c$, an operation covered by
Lemma~\ref{lem:closure}. The denominator has constant coefficient $1$,
so Theorem~\ref{thm:division} constructs an automaton for $a$.
If the numerator has $M$ kernel sequences, its polynomial states are
tuples $(P_1,\ldots,P_M)$ with $\deg P_j\le p^{e-1}$.
There are at most $p^{eM(p^{e-1}+1)}$ such tuples. Their transitions
are given by \eqref{eq:transition} with a single digit $r\in\{0,\ldots,p-1\}$.

For constant $b(k)=\beta$, the solution is $a(k)=c+k\beta$ in $R$,
which has period dividing $p^e$. With a general automatic $b$, the
same construction computes its partial sums modulo $p^e$.

\section{The constructive proof for nonnegative predecessors}\label{sec:nonnegative}

Throughout this section assume $P\subset\N^n\setminus\{0\}$.
Section~\ref{sec:negative} will supply the additional argument for
the full hypothesis \eqref{eq:geometry}.

\subsection{A product of equal-order cyclic groups}\label{sec:boundary-correction}

First take $H=R^s$, $R=\Z/p^e\Z$. A homomorphism $h:H^m\to H$ has
the form
\[
 h(z_1,\ldots,z_m)=\sum_{j=1}^m A_jz_j,
 \qquad A_j\in M_s(R),
\]
with column-vector conventions. Let $c_0$ be the automatic zero
extension of the boundary, as in Lemma~\ref{lem:boundary}.
Define the boundary correction $\gamma$ and the source array $r$ by
\begin{equation}\label{eq:boundary-source}
 \gamma(\x)=\one_C(\x)
      \left(c_0(\x)-\sum_j A_j(T_{\vct_j}c_0)(\x)\right),
 \qquad r(\x)=\gamma(\x)+\one_I(\x)b(\x).
\end{equation}
On $C$, the correction records the difference between the prescribed
value and the sum of its predecessor contributions, with values
outside $\N^n$ set to zero. On $I$, it vanishes. Its support can be
infinite, since a boundary face is infinite when $n\ge2$.
Every coordinate of $\gamma$ and $r$ is automatic by the closure lemmas.

If $\x\in C$ and $\x-\vct_j\in\N^n$, then
$\x-\vct_j\in C$. Indeed, some $x_i<d_i$, and
$x_i-(v_j)_i\le x_i<d_i$. This is precisely where the nonnegative
predecessor hypothesis enters the proof.

Extend $a$ by zero outside $\N^n$. The recurrence at interior points
and the preceding boundary observation give, at every point of $\N^n$,
\[
 a(\x)-\sum_j A_ja(\x-\vct_j)=r(\x).
\]
Let $S=\sum_{\x} a(\x)\X^{\x}$ and $R_0=\sum_{\x} r(\x)\X^{\x}$, as
column vectors of formal series. With
\[
 D(\X)=I_s-\sum_j A_j\X^{\vct_j},\qquad Q(\X)=\det D(\X),
\]
we have $DS=R_0$. Since every $\vct_j$ is nonnegative and nonzero,
$D(0)=I_s$ and $Q(0)=1$. 
By Lemmas~\ref{lem:adjugate} and~\ref{lem:seriesunit}, the
adjugate identity applies over $R[[\X]]$ and $Q$ is invertible. Thus
\begin{equation}\label{eq:solution}
 S=\frac{\adj(D)R_0}{Q}.
\end{equation}

Each coordinate $U_i$ of $\adj(D)R_0$ has automatic coefficients:
multiplication by a polynomial is a finite sum of translated zero
extensions with scalar coefficients. Apply
Theorem~\ref{thm:division} to every $U_i/Q$. Taking the tuple of the
coordinate outputs proves automaticity of $a$ and gives a terminating
construction.

\begin{corollary}[A computable bound]\label{cor:bound}
In this equal-exponent case, let $M_i$ count the states of an
automaton representing the sections of $U_i$, with transition
$\Lambda_{\rct}$ on input $\rct$. Let $q_j=\deg_{X_j}\det D$, and put
$V=\prod_j(p^{e-1}q_j+1)$. The vector solution has a
least-significant-digit automaton with at most
\[
 p^{eV\sum_{i=1}^s M_i}
\]
states. The quantities $M_i$ are computable from automata for the
boundary data and the additional sequence by the constructions above.
The bound applies before removing unreachable states or merging
equivalent states.
\end{corollary}
\begin{proof}

For coordinate $i$, Theorem~\ref{thm:division} gives a machine with
at most $p^{eVM_i}$ states. All these machines read the same digit
vectors from least to most significant. Their product, as defined
in Section~\ref{sec:product-automata}, has at most
\[
 \prod_{i=1}^s p^{eVM_i}=p^{eV\sum_{i=1}^s M_i}
\]
states and outputs the vector of solution coordinates.
\par
\end{proof}

\subsection{Different cyclic exponents}

Write $H=\prod_{i=1}^s\Z/p^{e_i}\Z$ and set $e=\max_i e_i$.
There is a coordinatewise quotient
$\pi:(\Z/p^e\Z)^s\to H$. Decompose $h$ into the sum of the
endomorphisms on its individual input factors. For each such
endomorphism use an integer matrix $A=(a_{ij})$ with
\begin{equation}\label{eq:divisibility}
 p^{e_i}\mid p^{e_j}a_{ij}.
\end{equation}
This condition is necessary and sufficient for multiplication by the
matrix to be well-defined on $H$. The same integer matrix modulo
$p^e$ defines an endomorphism of $(\Z/p^e\Z)^s$ commuting with $\pi$.
To see this, changing a representative in coordinate $j$ by
$p^{e_j}$ changes coordinate $i$ of its image by a multiple of
$p^{e_i}$, by \eqref{eq:divisibility}.

Choose one lift in $(\Z/p^e\Z)^s$ of every element of $H$. Lift $c$
and $b$ by this coding, which preserves automaticity and
compatibility on face intersections. The equal-exponent result gives
an automatic lifted solution $a'$. Applying $\pi$ to its recurrence
gives the original recurrence and boundary values. Uniqueness in
Lemma~\ref{lem:unique} gives $\pi a'=a$. A final output coding proves
Theorem~\ref{thm:abelian} in the nonnegative case.

\section{Allowing negative entries in one coordinate}\label{sec:negative}

With a negative predecessor coordinate, extending the prescribed boundary
by zero does not determine the whole boundary term in the generating
function equation. Some of that term depends on the unknown solution.
Condition~\eqref{eq:geometry} confines these unknown coefficients to
bounded degree in one variable. Algebraic division can then determine
them together with the solution.

For a field $K$, a formal series is \emph{algebraic} if it satisfies
a nonzero polynomial equation over the rational-function field in its
variables. We use the scalar algebraic division principle underlying
\cite[Theorem~13]{BMP} and its treatment in \cite{HK}, both in
characteristic zero.
The Christol--Salon theorem converts automaticity into algebraicity.
Algebraic Weierstrass division then determines the unknown boundary
series in the matrix argument.

\begin{theorem}[Christol--Salon]\label{thm:salon}
Let $K$ be a finite field of characteristic $p$.
A formal series $F\in K[[\X]]$ is algebraic over $K(\X)$ if and
only if its coefficient array is $p$-automatic.
\end{theorem}

For one variable, this is the theorem of Christol, Kamae,
Mend\`es France and Rauzy \cite[Theorem~1]{CKMR}.
Salon \cite{Salon} gives the multivariate extension. See also
\cite{DL} and the two-dimensional statement in
\cite[Theorem~14.4.1]{AS}.
Automaticity in base $p^f$ is equivalent to automaticity
in base $p$, because a finite automaton can group or ungroup blocks of
$f$ digits, allowing leading zero padding.

\begin{theorem}[Algebraic Weierstrass division]\label{thm:weierstrass}
Let $K$ be a field, write $Y=(Y_1,\ldots,Y_{n-1})$, and suppose
$D(Y,Z)$ is algebraic with $D(0,Z)=Z^L u(Z)$, where $u(0)\ne0$.
For every algebraic series $V(Y,Z)$ there are unique formal series
$Q$ and $R$ with
\[
 V=DQ+R,\qquad R\in K[[Y]][Z],\quad \deg_Z R<L.
\]
The quotient and all coefficients of $R$ are algebraic over the
corresponding rational-function fields. For $L=0$ the remainder is zero.
\end{theorem}

Formal Weierstrass division gives existence and uniqueness of $Q$
and $R$ among formal series. Algebraicity follows from \cite[Theorem~4.5]{Rond},
which holds over any field. In positive characteristic, the case of a
Weierstrass polynomial with inseparable factors is treated in
\cite[Lemma~4.3]{Rond}. We write $\operatorname{rem}_D V$
for the remainder. Formal uniqueness makes this operation $K[[Y]]$-linear: multiplying
a division identity by a series in $Y$ preserves the degree restriction
on the remainder, and adding division identities does the same.

\begin{lemma}[A matrix division step]\label{lem:matrixdivision}
Let $M(Y,Z)$ be an $s\times s$ polynomial matrix over $K$ satisfying
\[
 M(0,Z)=Z^\ell U(Z),\qquad U(0)=I_s.
\]
Suppose a vector of formal series $F$ and a vector
$R\in K[[Y]][Z]^s$ of degree less than $\ell$ satisfy
\begin{equation}\label{eq:matrixdivision}
 MF+R=V,
\end{equation}
where the entries of $V$ are algebraic. Then the entries of $F$ and
the coefficients of $R$ are algebraic.
\end{lemma}
\begin{proof}
For $\ell=0$, $M(0,0)=I_s$ and $F=M^{-1}V$ is algebraic.
Assume $\ell>0$. Put $D=\det M$ and $L=s\ell$. Then
$D(0,Z)=Z^L\det U(Z)$, so Theorem~\ref{thm:weierstrass} applies.
Write
\[
 R=\sum_{a=1}^s\sum_{j=0}^{\ell-1} r_{a,j}(Y)e_aZ^j,
\]
where $e_a$ are the standard coordinate vectors. By
Lemma~\ref{lem:adjugate}, multiplication by the adjugate replaces
the matrix coefficient $M$ by the scalar $D$. Multiply
\eqref{eq:matrixdivision} by $\adj M$ and take the scalar remainder
of each entry. This gives the finite linear system
\begin{equation}\label{eq:remaindersystem}
 \sum_{a,j} r_{a,j}(Y)
       \operatorname{rem}_D\bigl((\adj M)e_a Z^j\bigr)
   =\operatorname{rem}_D\bigl((\adj M)V\bigr).
\end{equation}
Equating coefficients of $Z^0,\ldots,Z^{L-1}$ in each component
gives $sL$ equations for $s\ell$ unknown series over $K[[Y]]$.
The $K[[Y]]$-linearity of the remainder justifies this step while
the unknown coefficients $r_{a,j}$ are still only formal series.
Every matrix entry and right-hand entry is an algebraic series in $Y$,
by Theorem~\ref{thm:weierstrass}.

This system has a square minor that is a unit at $Y=0$.
To prove this, specialization of a formal division identity at $Y=0$
shows that its remainder is truncation modulo $Z^L$. Also
\[
 (\adj M)(0,Z)=Z^{\ell(s-1)}\adj U(Z).
\]
If a vector polynomial $T(Z)$ of degree less than $\ell$ belongs to
the kernel of the specialized coefficient matrix, then
$Z^L$ divides $Z^{\ell(s-1)}\adj U(Z)T(Z)$. Since
$\adj U(0)=I_s$, Lemmas~\ref{lem:adjugate} and~\ref{lem:seriesunit}
give an inverse over $K[[Z]]$. Multiplication by that inverse shows
that $Z^\ell$ divides $T(Z)$. This forces $T=0$. The specialized columns
are therefore linearly independent over $K$, giving the asserted minor.

Solve the equations in this minor by Cramer's rule. Each $r_{a,j}$
is a fraction of algebraic series with a unit denominator, and is
therefore an algebraic power series. The assumed formal solution of
the full system agrees with this unique solution of the chosen equations.
Finally
\[
 F=\frac{(\adj M)(V-R)}{D}
\]
is algebraic componentwise. These fractions are power series because
they equal the original formal vector $F$. No regularity of the other
minors is needed.
\end{proof}

\subsection{The finite-field case}

Let $H=\F_p^s$, and write the rule as $h=\sum_j A_j$ on its separate
inputs. Let $c_0$ be the prescribed boundary extended by zero, and
put $f(u)=a(\boldsymbol d+u)$ for $u\in\N^n$, with $f=0$ outside
$\N^n$. Extend $a$ and $c_0$ similarly. Then everywhere on $\Z^n$,
\[
 a(x)=c_0(x)+f(x-\boldsymbol d).
\]
In particular, at every interior point the original recurrence is
equivalent to
\begin{equation}\label{eq:normalized}
 f(u)-\sum_j A_j f(u-\vct_j)=g(u)\quad(u\in\N^n),\qquad
 g(u)=b(\boldsymbol d+u)+\sum_j A_j c_0(\boldsymbol d+u-\vct_j).
\end{equation}
All terms defining $g$ are known and automatic, by the closure lemmas.

\phantomsection\label{sec:negative-boundary}
To see why the translation is needed, consider the boundary correction
that the original solution would require in a generating-function
equation on the whole orthant:
\[
 \gamma_a(x)=\one_C(x)
       \left(a(x)-\sum_j A_j a(x-\vct_j)\right).
\]
No recurrence is imposed on $C$, so this difference need not vanish.
For nonnegative predecessors, the observation in
Section~\ref{sec:boundary-correction} expresses it using $c_0$ alone.
With negative entries, a predecessor of a boundary point may lie in
$I$, where its value is still unknown. For example, take
$P=\{(1,-1),(0,1)\}$ and scalar coefficients $A_1=A_2=1$.
Here $\boldsymbol d=(1,1)$, and at the boundary point $(2,0)$,
\[
 \gamma_a(2,0)=c(2,0)-a(1,1),
\]
because the other predecessor $(2,-1)$ lies outside the orthant.
The value $a(1,1)$ is an interior value. Formula~\eqref{eq:normalized}
avoids putting such unknown values into the prescribed sequence $g$.
They remain part of the unknown series $F$ and the remainder $R$
in the equation below.

Name the exceptional coordinate $Z$, and the remaining coordinates
$Y$. This renaming does not change the causal order used to define the
array. In these variable names, write $\vct_j=(w_j,t_j)$ with
$w_j\in\N^{n-1}$ and $t_j\in\Z$. Set
\[
 \ell=\max(0,-t_1,\ldots,-t_m),\qquad
 M(Y,Z)=Z^\ell I_s-\sum_j A_j Y^{w_j}Z^{\ell+t_j}.
\]
Let $F,G$ be the generating vectors of $f,g$. Comparing coefficients
of $Y^uZ^{\ell+v}$, for $u\in\N^{n-1}$ and $v\ge0$, in
\eqref{eq:normalized} gives
\begin{equation}\label{eq:normalizedseries}
 MF+R=Z^\ell G,
 \qquad R\in\F_p[[Y]][Z]^s,\quad\deg_Z R<\ell.
\end{equation}
There are no remainders in the $Y$ directions because every $w_j$ is
coordinatewise nonnegative. If $w_j=0$, lexicographic positivity of
$\vct_j$ forces $t_j>0$. Consequently
\[
 M(0,Z)=Z^\ell\left(I_s-\sum_{w_j=0}A_j Z^{t_j}\right),
\]
and the matrix in parentheses has value $I_s$ at zero.

The entries of $G$ are algebraic by Theorem~\ref{thm:salon}.
Lemma~\ref{lem:matrixdivision} applied to \eqref{eq:normalizedseries}
makes the entries of $F$ algebraic. Theorem~\ref{thm:salon} now makes
$f$ automatic. Boundary assembly and translation make $a$ automatic.
This proves the elementary abelian case with the full geometric
hypothesis.

\subsection{Lifting to prime powers}

We next treat $H=(\Z/p^e\Z)^s$ by induction on $e$.
Reduce the solution modulo $p$. Its boundary data, affine term and matrices
satisfy the finite-field result, so $\bar a$ is automatic. Choose the
coordinate representatives $0,\ldots,p-1$ as a coding and let
$u$ be the resulting automatic lift to $(\Z/p^e\Z)^s$.
There is a unique array $z$ over $(\Z/p^{e-1}\Z)^s$ with
$a=u+pz$. On the interior it satisfies
\[
 z(x)=\sum_j (A_j\bmod p^{e-1})z(x-\vct_j)+\beta(x),
 \qquad
 \beta(x)=\frac{\sum_j A_j u(x-\vct_j)+b(x)-u(x)}p.
\]
The numerator lies in $p(\Z/p^e\Z)^s$ because $\bar a$ satisfies
the reduced recurrence. Division by $p$ here means the inverse of
the bijection from $(\Z/p^{e-1}\Z)^s$ to that ideal, which is a
well-defined coding. Extend $\beta$ by zero on the boundary.
It is automatic, as are the face values $z=(c-u)/p$.
The induction hypothesis therefore makes $z$ automatic, and then
$a=u+pz$ is automatic. The case $e=1$ was proved above.

For unequal cyclic exponents, apply exactly the matrix lift and output coding in the preceding section. That argument uses the
equal-exponent conclusion and uniqueness, and requires no further
geometric assumption. This completes the proof of
Theorem~\ref{thm:abelian} under \eqref{eq:geometry}.

The explicit state bound in Corollary~\ref{cor:bound} concerns the
nonnegative construction. The algebraic division proof above establishes
automaticity without a corresponding state bound.

\section{Proof for arbitrary finite \texorpdfstring{$p$}{p}-groups}\label{sec:groups}

We reduce a group recurrence to a quotient and a recurrence in a central
subgroup, following the approach of
\cite[Theorem~4 and Proposition~8(3)]{Pivato}. A subgroup is central
if its elements commute with the whole group, and it is fully invariant
if every endomorphism maps it into itself. We need both properties:
invariance gives a quotient rule, while centrality lets us collect the
subgroup contributions in the product without changing their order
relative to the other factors. An elementary abelian $p$-group is a
finite vector space over $\F_p$ under addition, so its recurrence is
covered by the abelian theorem.

\begin{lemma}\label{lem:subgroup}
Every nontrivial finite $p$-group $G$ has a nontrivial central elementary
abelian subgroup $N$ that is preserved by every endomorphism of $G$.
\end{lemma}
\begin{proof}
A finite $p$-group has nontrivial center by its class equation.
Induction on its order, applied to the quotient by the center, shows
it is nilpotent. Writing $[g,h]=ghg^{-1}h^{-1}$, define
$\gamma_1(G)=G$ and let $\gamma_{i+1}(G)$ be the subgroup generated
by the commutators $[u,g]$ with $u\in\gamma_i(G)$ and $g\in G$.
Nilpotence means that this series reaches the identity subgroup.
Let $\gamma_c(G)$ be the last nontrivial term of its
lower central series. This subgroup is central and is preserved by
every endomorphism, since each term
of the lower central series is generated by iterated commutators.
The elements of $\gamma_c(G)$ killed by the $p$th-power map form a
nontrivial subgroup $N$. It is elementary abelian because
$\gamma_c(G)$ is abelian. Both defining conditions are preserved by
endomorphisms, so $N$ has the required invariance.
\end{proof}

\begin{proof}[Proof of Theorem~\ref{thm:groups}]
Induct on $|G|$, the trivial group being immediate. Take $N$ as in
Lemma~\ref{lem:subgroup}. If $N=G$, the group is elementary abelian
and Theorem~\ref{thm:abelian} applies. Otherwise put $\overline G=G/N$.
Every $\alpha_j$ preserves $N$, so the ordered-product rule descends
to the same kind of rule on $\overline G$. Projecting the initial data
and $b$ gives automatic
data in $\overline G$.
 By induction, the projected solution
$q(\x)=a(\x)N$ is automatic.

Choose a set-theoretic section $\sigma:\overline G\to G$ with
$\sigma(1)=1$. Write uniquely
$a(\x)=\sigma(q(\x))z(\x)$ with $z(\x)\in N$. At interior points set
\[
 \beta(\x)=\sigma(q(\x))^{-1}
 \alpha_1\bigl(\sigma(q(\x-\vct_1))\bigr)\cdots
 \alpha_m\bigl(\sigma(q(\x-\vct_m))\bigr)b(\x).
\]

The quotient recurrence shows that $\beta(\x)\in N$. Extend $\beta$
by $1$ on $C$. It is automatic: the translated quotient arrays are
automatic by Lemma~\ref{lem:closure}, and the expression uses only
functions between finite sets followed by the automatic interior mask.

Since $N$ is central and each $\alpha_j(N)\subset N$, substitution
into the recurrence gives
\[
 z(\x)=\prod_{j=1}^m\alpha_j|_N(z(\x-\vct_j))\,\beta(\x)
 \qquad(\x\in I).
\]
On $C$, $z=\sigma(q)^{-1}c$ has automatic face restrictions. The last
display is an abelian recurrence on the elementary abelian group $N$
with an automatic affine term. Theorem~\ref{thm:abelian} makes $z$ automatic.
The coding $(q,z)\mapsto\sigma(q)z$ now makes $a$ automatic.

\end{proof}

\begin{remark}[Effectiveness and scope]
For a finite group specified by its multiplication table, the subgroups,
quotients, sections, and induced maps in this induction can be found by
enumerating subsets of $G$ and maps between the resulting
finite sets. For nonnegative predecessors, combining this with
the abelian construction gives an explicit effective proof. Its cost is
not analyzed here. For the larger geometry, the proof uses the algebraic
division step of Section~\ref{sec:negative}. Cross-commutation of the
images of different endomorphisms is unnecessary: only the centrality
of $N$ is used when collecting its factors. The unrestricted extension
fails already in its abelian special case.
\end{remark}

\subsection{Proof of the two-dimensional corollary}\label{sec:planar-proof}

\begin{proof}[Proof of Corollary~\ref{cor:plane}]
Write the predecessors as $\vct_j=(u_j,v_j)$.
Lexicographic positivity means that either $u_j>0$, or $u_j=0$
and $v_j>0$. Thus $u_j\ge0$ for all $j$, and a negative entry
can occur only in the second coordinate. This verifies
\eqref{eq:geometry} with the same exceptional coordinate $k=2$
for the entire predecessor set. There is no restriction on the
magnitudes or number of its negative second coordinates.

The recurrence is well-founded even when a predecessor has a larger
second coordinate than the point being computed. Choose a positive
integer $L$ so large that $Lu_j+v_j>0$ for each $j$ with $u_j>0$.
There are finitely many such inequalities. If $u_j=0$, the same
inequality holds because $v_j>0$. Each predecessor of $(m,n)$
therefore has smaller weight $Lm+n$. The weights are nonnegative
integers on $\N^2$. For each integer $N\ge0$, the set
$\{(m,n)\in\N^2:Lm+n\le N\}$ contains finitely many points,
because $m\le N/L$ and $n\le N$.

Put $d_1=\max(0,\max_j u_j)$ and
$d_2=\max(0,\max_j v_j)$. For $m\ge d_1$ and $n\ge d_2$,
all predecessors lie in $\N^2$. The remaining points form the
boundary $C=\{m<d_1\}\cup\{n<d_2\}$. Its face restrictions
are the automatic sequences assumed in the corollary, and their
compatibility makes the prescribed function well-defined.
The induction on these sublevel sets in
Lemma~\ref{lem:unique} determines a unique array with the specified data.

The sign condition and the automatic data now satisfy all hypotheses
of Theorem~\ref{thm:groups}, which makes the ordered-product solution
$p$-automatic. In the abelian case, the same conclusion also follows
directly from Theorem~\ref{thm:abelian}. If the rule is a homomorphism
$G^m\to G$, the factor decomposition in Corollary~\ref{cor:hom}
puts it in the ordered-product form.
\end{proof}

For example, the predecessor $(1,-2)$ requests the value at
$(m-1,n+2)$. With the additional predecessor $(0,1)$, the weight
$3m+n$ decreases by $1$ under either step. The automaticity conclusion
applies to this system for any of the group rules and data in the
corollary.

\section{Removing coordinates without prescribed boundary faces}\label{sec:time-boundary}

The restriction to one possible negative coordinate does not cover all
automatic predecessor systems. Another sufficient family allows
negative entries in several coordinates when the boundary is the
union of the slices $\{t=j\}$ for $0\le j<d$. Here $t$ denotes
the first coordinate, and the other $n-1$ coordinates range over
$\N^{n-1}$. There are $d$ slices, each infinite when $n\ge2$.
We first prove the closure property used to change coordinates.
The distinguished value $0$ in the lemma can be replaced by the
identity when the output alphabet is a group.

\begin{lemma}[Integer affine substitution]\label{lem:affine}
Let $u:\N^k\to A$ be $p$-automatic, where $A$ is finite with a
distinguished element $0$. Extend $u$ by zero to $\Z^k$.
For any integer matrix $B\in M_{k,d}(\Z)$ and $c\in\Z^k$, the
array $v(x)=u(Bx+c)$ on $\N^d$ is $p$-automatic.
\end{lemma}
\begin{proof}
For $e\ge0$ and $0\le r_i<p^e$, write
$Br+c=p^e q+s$ coordinatewise, with $0\le s_i<p^e$.
The coordinates of $q$ range over a finite set depending on $B,c$,
since $r_i/p^e$ belongs to $[0,1)$ and $|c_i|/p^e\le|c_i|$.
Let $u_{e,s}(y)=u(p^e y+s)$ for $y\in\N^k$, again extended by
zero outside $\N^k$. Then
\[
 v(p^e x+r)=u_{e,s}(Bx+q).
\]
The equality remains valid if a coordinate of $Bx+q$ is negative,
because $0\le s_i<p^e$ makes the corresponding coordinate of
$p^e(Bx+q)+s$ negative as well. There are finitely many choices of
$u_{e,s}$ by automaticity, and finitely many choices of $q$.
Thus $v$ has a finite $p$-kernel.
\end{proof}

\pagebreak[3]
\begin{theorem}[Boundary on prescribed coordinate slices]\label{thm:time-boundary}
Let $P$ be a finite nonempty set of vectors
\[
 \vct_j=(r_j,-s_j)\in\Z\times\Z^{n-1},
 \qquad r_j>0,\quad s_j\in\N^{n-1}.
\]
Put $d=\max_j r_j$. Let $G$ be a finite $p$-group,
let $\alpha_1,\ldots,\alpha_m$ be endomorphisms of $G$, and let
$b:\N\times\N^{n-1}\to G$ be $p$-automatic.
Prescribe $p$-automatic data on each slice $t=j$, $0\le j<d$, as functions
of the other coordinates. Then the unique solution of
\[
 a(t,x)=\alpha_1(a(t-r_1,x+s_1))\cdots
         \alpha_m(a(t-r_m,x+s_m))b(t,x)\quad(t\ge d)
\]
on $\N\times\N^{n-1}$ is $p$-automatic.
\end{theorem}
\begin{proof}
Choose integers $c_i\ge0$ such that $c_i r_j\ge (s_j)_i$ for
all $i,j$. Define the invertible integer linear map
\[
 T(t,x)=(t,y),\qquad y_i=x_i+c_i t.
\]
It sends each predecessor to
$\vct'_j=(r_j,c r_j-s_j)$, a nonnegative nonzero vector.
It sends the original domain to
\[
 W=\{(t,y)\in\N^n:y_i\ge c_i t\text{ for all }i\}.
\]
The boundary widths of the transformed nonnegative recurrence are
$d'_0=d$ and
\[
 d'_i=\max_j\bigl(c_i r_j-(s_j)_i\bigr)\le c_i d.
\]
If $t\ge d$, every point of $W$ lies in the transformed interior,
since $y_i\ge c_i t\ge d'_i$.

Let $c_0$ be the original boundary data extended by the identity to
all of $\N^n$. It is automatic by Lemma~\ref{lem:boundary}.
Extend it and the original additional sequence $b$ by the identity outside $\N^n$.
On the transformed boundary prescribe $c'=c_0\circ T^{-1}$, and on
the whole transformed orthant use the additional sequence $b'=b\circ T^{-1}$.
Both arrays are automatic by Lemma~\ref{lem:affine}, and the boundary
values are compatible because they are restrictions of one array.
The nonnegative case of Theorem~\ref{thm:groups} gives an automatic
solution $F$ for these transformed data.

Set $\widetilde a(t,x)=F(T(t,x))$. This array is automatic by
Lemma~\ref{lem:affine}. If $t<d$, its value is the prescribed
original boundary value, because $T(t,x)$ lies in the transformed
boundary and $T^{-1}T(t,x)=(t,x)$. If $t\ge d$, the recurrence
for $F$ applies at $T(t,x)$ and
\[
 T(t,x)-\vct'_j=T(t-r_j,x+s_j).
\]
These predecessors belong to $W$, since $t-r_j\ge0$ and
$x+s_j\in\N^{n-1}$. Substituting them into the transformed rule
gives the original ordered-product rule for $\widetilde a$, with additional sequence $b(t,x)$. Uniqueness gives $\widetilde a=a$.
\end{proof}

For example,
\[
 P=\{(1,-1,0),(1,0,-1)\}
\]
has negative entries in two coordinates. The map
$T(t,x,y)=(t,x+t,y+t)$ gives transformed predecessors
$\{(1,0,1),(1,1,0)\}$. Theorem~\ref{thm:time-boundary} applies
to any of the stated group rules and automatic data on $t=0$.
The argument uses the prescribed boundary slices as well as the transformed
predecessors. A coordinate change for a more general system would also
have to control its other boundary faces.

Related automaticity results for cellular automata appear in
\cite{GNW,RYCA,Nesme}.

\subsection{Coordinates without prescribed boundary faces}\label{sec:mixed-geometry}

A negative predecessor entry causes a boundary difficulty only when
the same coordinate also has a positive predecessor entry. If all its
entries are nonpositive, its boundary width is zero. We can then add a
positive weighted sum of the other coordinates without losing a
prescribed boundary face. This extends the argument for prescribed coordinate slices.

\begin{proof}[Proof of Theorem~\ref{thm:mixed}]
Put
\[
 A=\{i:d_i>0\},\qquad U=\{1,\ldots,n\}\setminus A.
\]
Thus $(v_j)_i\le0$ for $i\in U$. Project the predecessors onto
the coordinates in $A$, retaining their original order. The projected
vectors are nonzero and lexicographically positive: the first nonzero
coordinate of an original predecessor is positive and therefore belongs
to $A$. Lemma~\ref{lem:unique} supplies positive integer weights
$w_i$, $i\in A$, such that
\[
 h(x)=\sum_{i\in A}w_ix_i,\qquad \rho_j=h(\vct_j)>0.
\]
For each $i\in U$, choose an integer $c_i\ge0$ with
$c_i\rho_j+(v_j)_i\ge0$ for all $j$.

Define the invertible integer linear map $T$ by
\[
 (Tx)_i=x_i\quad(i\in A),\qquad
 (Tx)_i=x_i+c_i h(x)\quad(i\in U).
\]
The function $h$ depends only on the unchanged coordinates. Thus the
inverse subtracts $c_i h(x)$ in the coordinates in $U$.
The transformed predecessors $\vct'_j=T\vct_j$ are nonnegative
in $U$ and retain their entries in $A$. They remain lexicographically
positive, since their first original positive coordinate is unchanged
and earlier coordinates can only increase. In $A$, a negative entry
can occur only in the one coordinate permitted by the hypothesis.
Consequently $P'=\{\vct'_j\}$ satisfies \eqref{eq:geometry}.

We next check the transformed boundary. The image of the original
domain is
\[
 W=T(\N^n)=\{y\in\N^n:y_i\ge c_i h(y)\text{ for }i\in U\}.
\]
Its boundary widths satisfy $d'_i=d_i$ for $i\in A$ and
\[
 d'_i=\max_j\bigl((v_j)_i+c_i\rho_j\bigr)
       \le c_i\sum_{a\in A}w_a d_a\quad(i\in U).
\]
Indeed, $(v_j)_i\le0$ for $i\in U$ and
$(v_j)_a\le d_a$ for $a\in A$.
If $x\in I$, then $h(x)\ge\sum_{a\in A}w_a d_a$, so
$(Tx)_i\ge d'_i$ also in $U$. If $x\in C$, one of the unchanged
coordinates in $A$ is below its boundary width. Therefore
\[
 Tx\in I'\quad\Longleftrightarrow\quad x\in I,
 \qquad x\in\N^n.
\]

Let $c_0$ be the original boundary function extended by the identity
off $C$, and extend both $c_0$ and $b$ by the identity outside
$\N^n$. Prescribe $c_0\circ T^{-1}$ on the transformed boundary
and use the additional sequence $b\circ T^{-1}$. These data are automatic by
Lemmas~\ref{lem:boundary} and~\ref{lem:affine}. Their face restrictions
are compatible because they come from one array.
Theorem~\ref{thm:groups} gives an automatic transformed solution $F$.

The pullback $\widetilde a=F\circ T$ is automatic by
Lemma~\ref{lem:affine}. On $C$ it has the prescribed values.
For $x\in I$, the transformed recurrence applies at $Tx$, and
\[
 Tx-\vct'_j=T(x-\vct_j)\in W.
\]
Hence that recurrence becomes the original ordered-product recurrence
for $\widetilde a$, with additional sequence $b(x)$. Uniqueness gives
$\widetilde a=a$.
\end{proof}

For a three-dimensional example with two prescribed boundary directions,
take
\[
 P=\{(1,-1,-1),(0,1,-1)\},\qquad \boldsymbol d=(1,1,0).
\]
Both the second and third coordinates contain negative entries, but
only the second contains entries of both signs. The choice
$h(x,y,z)=2x+y$ and
\[
 T(x,y,z)=(x,y,z+2x+y)
\]
gives transformed predecessors $\{(1,-1,0),(0,1,0)\}$.
Theorem~\ref{thm:mixed} allows arbitrary compatible automatic data on
$x=0$ and $y=0$, together with an automatic additional sequence. The earlier
family with boundary on prescribed slices has no coordinate with both signs and is a
special case of this theorem.

If no coordinate contains both signs, all transformed predecessors
are nonnegative. In that case the reduction combines with the direct
construction in Section~\ref{sec:nonnegative}, without algebraic division.
The integer affine substitution itself also admits an explicit machine,
as follows.

\begin{proposition}[A bound for integer affine substitution]\label{prop:affine-bound}
Suppose $u:\N^k\to A$ is given by a least-significant-digit automaton
with $M$ states, and extend it by a distinguished value outside
$\N^k$. For $v(x)=u(Bx+c)$ as in Lemma~\ref{lem:affine}, put
\[
 L_i=\min\left(c_i,\sum_{j:B_{ij}<0}B_{ij}\right),\qquad
 U_i=\max\left(c_i,\sum_{j:B_{ij}>0}B_{ij}\right).
\]
There is an effectively constructed automaton for $v$ with at most
$M\prod_{i=1}^k(U_i-L_i+1)$ states.
\end{proposition}
\begin{proof}
First discard unreachable input states, leaving at most $M$ states.
For each remaining state $q$, let $u_q$ be the array computed by
starting the input automaton in $q$. These arrays inherit padding
invariance, as in the translation construction in
Lemma~\ref{lem:closure}.
A state is a pair $(q,z)$, where $q$ is an input automaton state
and $L_i\le z_i\le U_i$. It represents the function
$x\mapsto u_q(Bx+z)$, with the same extension outside $\N^k$.
The initial pair is $(q_0,c)$. On an input digit vector $r$, divide
componentwise
\[
 Br+z=pz'+s,\qquad 0\le s_i<p.
\]
The successor is $(\delta(q,s),z')$. The represented function is the
required section, including at negative arguments, by the calculation
in Lemma~\ref{lem:affine}.
The interval for $z'_i$ is invariant because
\[
 pL_i\le (Br+z)_i\le pU_i,
\]
using $0\le r_j<p$, $L_i\le\sum_{B_{ij}<0}B_{ij}$ and
$U_i\ge\sum_{B_{ij}>0}B_{ij}$. Taking the integer quotient gives
$L_i\le z'_i\le U_i$.
The output is $u_q(z)$ if $z\in\N^k$, and the distinguished
value otherwise. Evaluating $u_q$ at the fixed tuple $z$ uses its
given automaton. There are at most $M$ choices of $q$ and
$\prod_{i=1}^k(U_i-L_i+1)$ choices of $z$, giving at most
$M\prod_{i=1}^k(U_i-L_i+1)$ states.
\end{proof}

\section{A three-dimensional obstruction}\label{sec:limit}
\label{sec:counterexample}

With negative entries in two coordinate directions, automatic boundary
data need not give an automatic solution. The construction uses the
recurrence of \cite[Section~4.5, Example~6]{BMP}. Its
generating-function identity is the Gabber--Kashiwara division of $XY$
by $(X-Y^2)(Y-X^2)$ in \cite[Example~5.2]{ACH}. An additional first coordinate makes
its predecessor vectors lexicographically positive and allows prescribed
boundary data with a single nonzero value.

\begin{theorem}\label{thm:counterexample}
For every prime $p$, there is a recurrence over the additive group
$\F_p$, with lexicographically positive predecessors, a homomorphism
rule and a single nonzero boundary value, whose solution on $\N^3$
is not $p$-automatic.
\end{theorem}

\subsection{The recurrence and its initial conditions}

Use coordinates $(t,i,j)$, in that order, and set
\[
 P=\{(1,2,-1),(1,-1,2),(2,1,1)\}.
\]
Every vector has positive first coordinate. The prescribed boundary is
\[
 C=\{(t,i,j)\in\N^3:t<2\text{ or }i<2\text{ or }j<2\}.
\]
Prescribe
\[
 c(t,i,j)=\begin{cases}1,&(t,i,j)=(0,2,2),\\0,&\text{otherwise},\end{cases}
 \qquad (t,i,j)\in C.
\]
For $t,i,j\ge2$ define, in $\mathbb F_3$,
\begin{equation}\label{eq:counterrec}
 a(t,i,j)=a(t-1,i-2,j+1)+a(t-1,i+1,j-2)-a(t-2,i-1,j-1).
\end{equation}
The rule $(u,v,w)\mapsto u+v-w$ is a group homomorphism
$(\mathbb F_3)^3\to\mathbb F_3$. Every specified boundary slice is zero
except the slice $t=0$, which has a single nonzero point. Such slices are
automatic in every base: a finite automaton can check equality with a
fixed coordinate tuple. The values agree on all intersections. The
recurrence has a unique solution by induction on $t$, since each right
hand side uses points with a smaller first coordinate.

\subsection{The formula for the section}\label{sec:section-formula}

\begin{proposition}\label{prop:section-formula}
The solution of \eqref{eq:counterrec} satisfies
\begin{equation}\label{eq:section}
 a(n,n+2,2)=
 \begin{cases}(-1)^k\bmod3,&n=3\cdot2^k-3,\ k\ge0,\\
 0,&\text{otherwise}.
 \end{cases}
\end{equation}
\end{proposition}
\begin{proof}
First work over the integers. Define $q_{ij}=0$ when either index is
negative, and for $i,j\ge0$ set
\begin{equation}\label{eq:qrec}
 q_{ij}=\mathbf1_{(i,j)=(0,0)}
       +q_{i-2,j+1}+q_{i+1,j-2}-q_{i-1,j-1}.
\end{equation}
Every term with nonnegative indices on the right has smaller total degree.
This defines $q$ uniquely. In particular $q_{00}=1$ and $q_{10}=q_{01}=0$.
Let
\[
 Q(X,Y)=\sum_{i,j\ge0}q_{ij}X^iY^j,\qquad
 g(X,Y)=(X-Y^2)(Y-X^2).
\]
The recurrence and initial value are symmetric in $i,j$, so
$Q(X,Y)=Q(Y,X)$. For $u,v\ge1$, coefficient extraction gives
\[
 [X^uY^v]gQ
 =q_{u-1,v-1}-q_{u-3,v}-q_{u,v-3}+q_{u-2,v-2}
 =\mathbf1_{u=v=1}.
\]
On the axes, $g(X,0)Q(X,0)=-X^3Q(X,0)$, and likewise for $X=0$.
Consequently, writing $r(X)=X^3Q(X,0)$, we obtain the exact formal identity
\begin{equation}\label{eq:gq}
 g(X,Y)Q(X,Y)=XY-r(X)-r(Y).
\end{equation}
Substituting $Y=X^2$ is valid in formal power series, since a coefficient
of $X^n$ receives contributions from only finitely many pairs $(i,j)$.
It yields
\[
 r(X)+r(X^2)=X^3.
\]
For any positive integer $N$, repeated substitution gives
\[
 r(X)=\sum_{k=0}^{N-1}(-1)^kX^{3\cdot2^k}
          +(-1)^N r(X^{2^N}).
\]
The last summand has order at least $3\cdot2^N$ because
$r(X)=X^3Q(X,0)$. Thus a fixed coefficient is unaffected by that
summand when $N$ is sufficiently large, proving
\begin{equation}\label{eq:lacunary}
 r(X)=\sum_{k\ge0}(-1)^kX^{3\cdot2^k},\qquad
 Q(X,0)=\sum_{k\ge0}(-1)^kX^{3\cdot2^k-3}.
\end{equation}

The solution of \eqref{eq:counterrec} is
\begin{equation}\label{eq:embedding}
 a(t,i,j)=
 \begin{cases}
 q_{i-2,j-2}\bmod3,&i,j\ge2\text{ and }t=i+j-4,\\
 0,&\text{otherwise}.
 \end{cases}
\end{equation}
We verify both the recurrence and its prescribed values. Put
$\Delta(t,i,j)=t-i-j+4$. For the three predecessor vectors,
\[
 1-2-(-1)=0,\qquad 1-(-1)-2=0,\qquad 2-1-1=0.
\]
Hence $\Delta$ is preserved by each predecessor step. At an interior
point with $\Delta\ne0$, all terms in \eqref{eq:counterrec} vanish
under the proposed formula. At an interior point with $\Delta=0$,
put $u=i-2$ and $v=j-2$. The three predecessor values become
$q_{u-2,v+1}$, $q_{u+1,v-2}$ and $q_{u-1,v-1}$, reduced modulo $3$.
A negative index gives a boundary value zero. Since $t=u+v\ge2$,
the seed term of \eqref{eq:qrec} is zero, and that equation proves
\eqref{eq:counterrec} at the point.

The boundary cases also agree. If $i<2$ or $j<2$, the formula is zero.
If $t=0$, the conditions $i,j\ge2$ and $i+j=4$ force $i=j=2$,
where the value is $q_{00}=1$. If $t=1$, the possible nonzero entries
would be $q_{10}$ and $q_{01}$, both zero. Thus the formula has the
prescribed values on all boundary faces. Uniqueness of the recurrence
proves \eqref{eq:embedding}.

At $(t,i,j)=(n,n+2,2)$, the equality $t=i+j-4$ holds and the
embedded coefficient is $q_{n,0}$. Extracting that coefficient from
\eqref{eq:lacunary} gives \eqref{eq:section}.
\end{proof}
The nonzero positions of this section begin with
$0,3,9,21,45,93,189,381$.

\subsection{Why this section cannot be automatic in base three}

A subset of $\N$ is \emph{recognizable in base $k$} if a finite
automaton accepts its canonical base-$k$ representations. We show
that the support in \eqref{eq:section} is not recognizable in base $3$,
using a repeated cycle in an accepting path.

\begin{lemma}\label{lem:sparse}
Let $b,k\ge2$ be multiplicatively independent, meaning that
$b^u\ne k^v$ for positive integers $u,v$. For any positive integer $C$, the
set $\{C(b^j-1):j\ge0\}$ is not recognizable in base $k$.
\end{lemma}
\begin{proof}
Suppose a finite automaton accepts the canonical base-$k$
representations. There are infinitely many such words, so an accepting
path for a sufficiently long word
contains a cycle after its first digit. There are fixed words $u,v,w$,
with $u$ starting in a nonzero digit and $d=|v|>0$, such that every
$uv^\ell w$, $\ell\ge0$, is accepted. Let $N_\ell$ be its numerical
value. Summing the geometric progression for the repeated block gives
\[
 N_\ell=A k^{d\ell}+B
\]
for fixed rational $A>0,B$. Hence
\[
 \frac{N_{\ell+1}+C}{N_\ell+C}\longrightarrow k^d.
\]
Since each $N_\ell+C$ is $Cb^{j_\ell}$, each ratio is an integral
power of $b$. In a compact subinterval of $(0,\infty)$ there are only
finitely many integral powers of $b$. Convergence therefore forces the
ratios eventually to equal their limit. Thus $k^d=b^h$ for a positive
integer $h$, a contradiction.
\end{proof}

If $a$ were $3$-automatic, Lemma~\ref{lem:affine} applied to the
integer affine map $n\mapsto(n,n+2,2)$ would make the section
$s(n)=a(n,n+2,2)$ $3$-automatic. Coding its nonzero values as $1$ gives a
base-three recognizer for its support. Formula \eqref{eq:section} and
Lemma~\ref{lem:sparse}, with $b=2,k=3,C=3$, give a contradiction.
This proves the assertion of Theorem~\ref{thm:counterexample} for $p=3$.

\subsection{The obstruction occurs for every prime}

For completeness, replace the integer $2$ controlling the dilation by
an integer $b\ge2$, and take
\[
 P_b=\{(b-1,b,-1),(b-1,-1,b),(2b-2,b-1,b-1)\}.
\]
The initial domain is $t<2b-2$ or $i<b$ or $j<b$. Set the single initial
value at $(0,b,b)$ to $1$ and all other initial values to $0$. The rule
remains $(u,v,w)\mapsto u+v-w$. Define the integer array
\[
 q_{ij}=\mathbf1_{(i,j)=(0,0)}+q_{i-b,j+1}
       +q_{i+1,j-b}-q_{i-b+1,j-b+1},
\]
with negative indices zero. Its terms are determined by smaller total
degree. For $0<i+j<2b-2$, all values vanish: an induction leaves only
possible predecessors $(0,0)$, and reaching that point with either of
the first two shifts would require a negative coordinate, while the
third shift requires total degree $2b-2$.

The preceding coefficient calculation, now with
$g_b=(X-Y^b)(Y-X^b)$, gives
\[
 g_bQ=XY-r_b(X)-r_b(Y),\quad r_b(X)=X^{b+1}Q(X,0),\quad
 r_b(X)+r_b(X^b)=X^{b+1}.
\]
It follows that
\[
 q_{n,0}=\begin{cases}(-1)^j,&n=(b+1)(b^j-1),\\0,&\text{otherwise}.
 \end{cases}
\]
The embedding is $a(t,i,j)=q_{i-b,j-b}$ on $t=i+j-2b$, zero elsewhere,
reduced modulo $p$. Its initial layers agree by the vanishing just
proved. Its section $a(n,n+b,b)$ has support
$\{(b+1)(b^j-1):j\ge0\}$. If $a$ were $p$-automatic,
Lemma~\ref{lem:affine} applied to $n\mapsto(n,n+b,b)$ would make this
section $p$-automatic. Lemma~\ref{lem:sparse} applies whenever $b$ and $p$ are
multiplicatively independent. Choosing $b=2$ for odd $p$, and $b=3$ for
$p=2$, proves failure of the unrestricted assertion for every prime,
already with the cyclic group $\mathbb F_p$.

The example also gives failure in every dimension $n>3$.
Append $n-3$ zero coordinates to each predecessor and prescribe the
same boundary values independently of the added indices. Uniqueness
makes the new solution independent of those indices as well.
If it were $p$-automatic, feeding zero digits in the added coordinates
would give an automaton for the original three-dimensional array,
which is impossible.

\subsection{Two different jumps and a fixed obstruction for every prime}
\label{sec:asymmetric-obstruction}

The preceding symmetric example changes its jump size when the prime
changes. Allowing two different jump sizes gives a single recurrence
that fails in every prime characteristic. It also identifies when the
same factorized construction becomes automatic. The factorization is
an adaptation of \cite[Section~4.5, Example~6]{BMP}, not a separate
source of the sparse-series method.

We first record the support test used in the nonautomatic case.
Here the support of an array is the set of indices where its value
is nonzero.

\begin{lemma}[A finite union of geometric progressions]\label{lem:geometric-union}
Let $L,k\ge2$ be multiplicatively independent, let
$A_1,\ldots,A_s$ be positive integers, and let $C\in\Z$.
No infinite subset of
\[
 \left\{A_iL^j-C:1\le i\le s,\ j\ge0\right\}\cap\N
\]
is recognizable in base $k$.
\end{lemma}
\begin{proof}
If such a subset were recognizable, the accepting-cycle argument in
Lemma~\ref{lem:sparse} would give accepted numbers
$N_\ell=a k^{d\ell}+b$, where $a>0$, $d>0$ and $a,b$ are fixed
rationals. Thus
\[
 \frac{N_{\ell+1}+C}{N_\ell+C}\longrightarrow k^d.
\]
For sufficiently large $\ell$, each ratio belongs to
$\{(A_i/A_j)L^h:1\le i,j\le s,\ h\in\Z\}$.
This set has finite intersection with any compact subinterval of
$(0,\infty)$. Convergence forces the ratios to equal $k^d$
eventually. Among the finitely many possible multipliers $A_i$, one
occurs at two later indices $\ell$ and $\ell+r$, where $r>0$.
Dividing those two terms gives $k^{dr}=L^h$ for an integer $h>0$,
contrary to multiplicative independence.
\end{proof}

For the positive case over prime powers, we need division by a
polynomial whose constant coefficient may vanish. The existence of
the formal quotient is an assumption in this lemma.

\begin{lemma}[Exact polynomial division]\label{lem:exact-division}
Let $R=\Z/p^e\Z$, and let $D\in R[\X]$ have nonzero reduction
modulo $p$. If $F\in R[[\X]]$ and the coefficients of $DF$ are
$p$-automatic, then the coefficients of $F$ are $p$-automatic.
\end{lemma}
\begin{proof}
For $e=1$, Theorem~\ref{thm:salon} makes $DF$ algebraic, so
$F=(DF)/D$ is algebraic over $\F_p(\X)$. It is already a formal
series by assumption, and Theorem~\ref{thm:salon} applies again.
For $e>1$, reduction modulo $p$ therefore gives an automatic series
$\overline F$. Lift its coefficients to $\{0,\ldots,p-1\}$ to
obtain an automatic series $U$ over $R$. Write $F=U+pZ$, with $Z$
over $\Z/p^{e-1}\Z$. The coefficients of
\[
 DZ=\frac{DF-DU}{p}
\]
are automatic. Multiplication by $D$ uses a sum of translated
coefficient arrays, one for each monomial of $D$. Subtraction and
division by $p$ on the ideal $pR$ are codings between finite alphabets.
Lemma~\ref{lem:closure} therefore preserves automaticity in these
operations. The reduced polynomial still has nonzero
reduction modulo $p$, so induction makes $Z$ automatic.
The coding $F=U+pZ$ finishes the proof.
\end{proof}

\begin{theorem}[The two-parameter family]\label{thm:asymmetric}
Let $b,c\ge2$, let $p$ be prime, and let $e\ge1$.
Use predecessors
\[
 P_{b,c}=\{(b-1,b,-1),(c-1,-1,c),(b+c-2,b-1,c-1)\}
\]
and the additive rule $(u,v,w)\mapsto u+v-w$ over $\Z/p^e\Z$.
On the boundary
$\{t<b+c-2\}\cup\{i<b\}\cup\{j<c\}$, prescribe $1$ at
$(0,b,c)$ and $0$ elsewhere. The unique solution is $p$-automatic
if and only if $bc$ is a power of $p$.
\end{theorem}
\begin{proof}
Work first over the integers. Put $q_{ij}=0$ at negative indices and
define
\[
 q_{ij}=\mathbf1_{i=j=0}+q_{i-b,j+1}+q_{i+1,j-c}
                            -q_{i-b+1,j-c+1}.
\]
All right-hand indices in the nonnegative quadrant have smaller
total degree. For $0<i+j<b+c-2$, induction gives $q_{ij}=0$.
Indeed, reaching the seed through either of the first two shifts
would require a negative coordinate, and reaching it through the
third requires $(i,j)=(b-1,c-1)$.

Set $Q(X,Y)=\sum q_{ij}X^iY^j$ and
\[
 g=(X-Y^c)(Y-X^b),\qquad
 R(X)=X^{b+1}Q(X,0),\qquad S(Y)=Y^{c+1}Q(0,Y).
\]
Coefficient extraction as in \eqref{eq:gq} gives
\[
 gQ=XY-R(X)-S(Y).
\]
Substitute $Y=X^b$ and $X=Y^c$ in this identity. Both substitutions
are valid formally. They give
\[
 R(X)+S(X^b)=X^{b+1},\qquad
 R(Y^c)+S(Y)=Y^{c+1}.
\]
Eliminating $S$ yields, with $L=bc$,
\[
 R(X)-R(X^L)=X^{b+1}-X^{b(c+1)}.
\]
Iteration is justified coefficient by coefficient because the order
of the remaining term tends to infinity. Therefore
\begin{equation}\label{eq:asymmetric-series}
 R(X)=\sum_{r\ge0}\left(X^{(b+1)L^r}-X^{b(c+1)L^r}\right).
\end{equation}
The corresponding formula for $S$ exchanges $b$ and $c$.
The two exponent sets in \eqref{eq:asymmetric-series} are disjoint:
\[
 1<\frac{b(c+1)}{b+1}<bc=L.
\]
Thus their coefficients remain nonzero modulo every $p^e$.

The three-dimensional solution is
\[
 a(t,i,j)=
 \begin{cases}q_{i-b,j-c}\bmod p^e,&i\ge b,\ j\ge c,
                       \ t=i+j-b-c,\\0,&\text{otherwise}.
 \end{cases}
\]
For each predecessor, its first coordinate is the sum of its other
two coordinates. The equation $t=i+j-b-c$ is consequently preserved
under predecessor steps, and the recurrence on that plane becomes
the displayed recurrence for $q$. The proved vanishing for
$0<i+j<b+c-2$ supplies the prescribed slices $0\le t<b+c-2$. The other
boundary faces are zero. Uniqueness therefore proves this embedding.

If $L$ is not a power of $p$, it is multiplicatively independent of
$p$. The section $a(n,n+b,c)$ has support
\begin{equation}\label{eq:asymmetric-support}
 \{(b+1)L^r-(b+1):r\ge0\}
 \ \cup\ \{b(c+1)L^r-(b+1):r\ge0\}.
\end{equation}
Lemma~\ref{lem:geometric-union} says that this support is not
recognizable in base $p$. If the full array were automatic, affine
substitution and a nonzero-value coding would make the support
recognizable. This contradiction proves nonautomaticity.

Suppose instead that $L=p^h$. The base-$p$ representations of
$Ap^{hr}$, for fixed $A>0$, are the representation of $A$ followed
by $hr$ zeros. A finite automaton recognizes them by checking the
fixed prefix and counting the appended zeros modulo $h$.
Formula~\eqref{eq:asymmetric-series} therefore gives automatic
coefficient arrays for $R$ and $S$ modulo $p^e$. So $XY-R-S$ is
automatic. The coefficient of $XY$ in $g$ is $1$, hence $g$ has
nonzero reduction modulo $p$. Lemma~\ref{lem:exact-division} makes
the coefficients of $Q$ automatic.
Finally, integer affine substitution and the automatic test
$t-i-j+b+c=0$ make the embedded array automatic. The equality test
is obtained by applying Lemma~\ref{lem:affine} to the indicator of
$\{0\}$, extended by zero to the integers.
\end{proof}

\begin{corollary}[One fixed predecessor set for all primes]
\label{cor:uniform-obstruction}
The fixed set
\[
 P=\{(1,2,-1),(2,-1,3),(3,1,2)\}
\]
with the rule $(u,v,w)\mapsto u+v-w$ and the single boundary seed
at $(0,2,3)$ produces a non-$p$-automatic array over $\Z/p^e\Z$
for every prime $p$ and every $e\ge1$.
\end{corollary}
\begin{proof}
Take $b=2,c=3$ in Theorem~\ref{thm:asymmetric}. The integer $6$
is not a power of any prime.
\end{proof}

For this set, the section $a(n,n+2,3)$ is $1$ at
$n=3\cdot6^r-3$, is $-1$ at $n=8\cdot6^r-3$, and is zero
elsewhere, with values reduced modulo $p^e$.
Both the second and third predecessor coordinates contain entries of both signs.
Together with Theorem~\ref{thm:mixed}, this shows that the threshold
of at most one such coordinate is sharp as a sufficient condition
depending only on their number. It does not classify individual
systems with two or more coordinates of both signs.

\section{Verifying and constructing automata}\label{sec:certified-construction}

For systems treated by algebraic power-series division, the
proof establishes automaticity without directly giving a transition
table. We can instead compute candidate
kernel states from finite samples and verify the recurrence for every
index. For one-dimensional sequences, kernel inference from finite
prefixes is established, including the possibility of false
identifications and tight bounds on the prefix length
\cite[Section~2.3 and Theorem~2.3]{Sutner}.
The verification below uses the specified recurrence to decide whether
a candidate computes its unique solution.

\begin{proposition}[Finite verification of a candidate]\label{prop:verify-candidate}
Given automata for the boundary data, the additional sequence, and a candidate array $u$,
one can decide whether $u$ satisfies the prescribed recurrence at all
indices. The procedure either certifies $u=a$ or returns an index where
the boundary condition or recurrence fails.
\end{proposition}
\begin{proof}
Construct automata for the interior mask and the assembled boundary
function using Lemmas~\ref{lem:closure} and~\ref{lem:boundary}.
For each predecessor, Proposition~\ref{prop:affine-bound}, with
$B=I$ and $c=-\vct_j$, constructs the shifted array
$u(x-\vct_j)$, extended by the neutral element at negative indices.

Use Lemma~\ref{lem:kernel} to put all the machines in the
least-significant-digit convention, and form their product as in
Section~\ref{sec:product-automata}. At a product state, let $z$ be
the candidate output, $z_j$ the output for predecessor $j$, $\beta$
the additional-sequence output, $\gamma$ the assembled boundary
output, and $\iota\in\{0,1\}$ the interior mask. For the ordered
group rule, give this state the Boolean output
\[
 \kappa(z,z_1,\ldots,z_m,\beta,\gamma,\iota)=
 \begin{cases}
 1,&\iota=0\text{ and }z=\gamma,\\
 1,&\iota=1\text{ and }
       z=\alpha_1(z_1)\cdots\alpha_m(z_m)\beta,\\
 0,&\text{otherwise}.
 \end{cases}
\]
This is a coding on a finite output alphabet. An arbitrary specified
local rule on a finite alphabet can be checked by substituting that
rule in the second comparison.

Enumerate the reachable product states by the construction in
Section~\ref{sec:product-automata}, recording a word reaching each
new state. If a state outputs $0$, its word
$\rct_0\cdots\rct_{\ell-1}$ gives the failing index
$x=\sum_{j<\ell}p^j\rct_j$. The empty word checks $x=0$.
If every reachable state outputs $1$, every index satisfies its
boundary condition or its interior recurrence, as appropriate.
Lemma~\ref{lem:unique} then gives $u=a$.
\par
\end{proof}

\begin{theorem}[A terminating search when the solution is automatic]
\label{thm:verified-search}
Given a recurrence with a finite nonempty set of lexicographically
positive predecessor offsets, a local rule on a finite output alphabet,
and automatic compatible data, there is an algorithm
that constructs and verifies an automaton whenever the unique solution
is $p$-automatic. In particular, it terminates under the hypotheses of
Theorem~\ref{thm:mixed}.
\end{theorem}
\begin{proof}
For a positive integer $B$, compare sections of the solution on the
finite set $\{0,\ldots,B-1\}^n$. A section is represented by a pair
$(e,r)$ and means the function $x\mapsto a(p^e x+r)$.
Each requested value can be computed from the recurrence and data by
induction on the positive weight in Lemma~\ref{lem:unique}.

Starting from $a$, process all digit sections of each retained
representative. Retain a new representative only if its vector of
sample values differs from those already retained. Otherwise use the
existing representative with that vector. There are at most
$|G|^{B^n}$ sample vectors when the output alphabet is $G$, so this
process terminates for each $B$. It defines a finite candidate machine,
with a state's output its representative's value at zero.
The sample contains zero, so the output is unchanged by zero padding.

Apply Proposition~\ref{prop:verify-candidate}. Return the candidate
if it passes. Otherwise increase $B$ and repeat.
If $a$ is automatic, its kernel is finite. For each pair of distinct
kernel arrays, choose one index where they differ. All these finitely
many indices lie in some cube $\{0,\ldots,B_0-1\}^n$.
For $B\ge B_0$, comparison on that cube identifies kernel arrays
correctly. The closure procedure then constructs the actual kernel
machine, which passes verification. An earlier candidate that passes
is also correct by uniqueness.
\end{proof}

The theorem does not give a stopping criterion for a nonautomatic
solution. It also gives no bound on the number of computational steps in terms of the
recurrence's input size. The coefficient samples propose a machine,
and checking all reachable product states supplies its certificate.

\subsection{Two certified binary examples}\label{sec:binary-certificates}

The first example has predecessors
$\{(1,0),(1,2),(2,-2)\}$ and satisfies
\begin{equation}\label{eq:certified-negative}
 a(m,n)=a(m-1,n)+a(m-1,n-2)+a(m-2,n+2)+1\pmod2
 \quad(m,n\ge2),
\end{equation}
with $a(m,n)=1$ when $m<2$ or $n<2$. The constant term matters:
if it is omitted, the same boundary data give the constant array $1$,
because $1+1+1=1$ in $\F_2$.
For \eqref{eq:certified-negative}, sample cubes of side $1,2,3$
give incorrect candidates. Verification returns failing indices
$(2,2)$, $(4,4)$ and $(52,16)$, respectively. Side $4$ gives a
candidate with $53$ states, which passes the verification.

The second example has
\[
 P=\{(0,4),(1,3),(3,1),(4,0),(4,4)\}
\]
and rule $a(x)=\sum_{v\in P}a(x-v)+1$ over $\F_2$, with
unit boundary data on $m<4$ or $n<4$. For this example, side $6$
gives a verified candidate with $254$ states.

Table~\ref{tab:certified-binary} gives state counts and block-rule
counts for these two specified infinite arrays. The digit automata
read the coordinate digits together, with zero padding allowed.
A numerical block is a $16\times16$ array at an origin
$(16m,16n)$. The numerical substitution replaces each block by the
corresponding $32\times32$ array. A colored substitution may use
different labels for the same numerical block, and distinct labels
may share a display color.

\begin{table}[!htbp]
\centering
\begin{tabular}{lrrrr}
\toprule
Example & Least digits first & Most digits first & Blocks & Labels\\
\midrule
Negative predecessor & $53$ & $107$ & $44$ & $44$\\
Five nonnegative predecessors & $254$ & $710$ & $204$ & $212$\\
\bottomrule
\end{tabular}
\caption{Counts for the two binary recurrences in
Section~\ref{sec:binary-certificates}. The first two numerical columns
give minimal automaton state counts for the indicated reading directions,
using binary digit pairs and allowing zero padding. The Blocks column
counts distinct $16\times16$ numerical blocks at origins $(16m,16n)$.
The Labels column counts states of the reconstructed colored substitution.
Both substitutions have scale factor $2$.}
\label{tab:certified-binary}
\end{table}

The verification programs, input data and computational certificates
are available in the accompanying repository.\footnote{\label{fn:repository}GitHub Repository: \url{https://github.com/bogdan27182/paper-ho}.}
All file paths in this paper are relative to the repository root. The
commands in \path{README.txt} reproduce the reported computations.

The script \path{research/certify_original_examples.py} produces the
full transition tables and certificates. The recurrence checks explore
$686$ reachable product states for the first example and $552\,493$
for the second. At each reachable state, the boundary or recurrence
comparison holds. The minimization procedure in
Section~\ref{sec:product-automata}, which refines states by their
outputs and successor classes, certifies the minimal state counts
in the table.

The block certificate examines every reachable state of each
most-significant-digit machine. At such a state, four further digit
pairs give its $16\times16$ block, and five give its $32\times32$
expansion. The certificate checks that equal smaller blocks always
have equal expansions, that all four child blocks belong to the
same collection of numerical blocks, and that the initial block is its own upper
left child. These identities make substitution iteration agree with
the array on squares of side $16\cdot2^r$ for every $r\ge0$.
The colored-state transitions decode to these numerical rules.

For reproduction, the two numerical images and the two colored images
distributed with the examples each have $512^2=262\,144$ pixels.
The numerical images match the respective recurrences at all displayed
positions. Every reachable most-significant-digit state occurs among
the displayed positions in each colored image, and its assigned color
is consistent across all occurrences. The resulting colored machines
decode to the numerical block rules at all indices by the preceding
finite checks. The $44$ and $212$ labels use $39$ and $151$ distinct
RGB colors, respectively, because the palettes repeat some colors.

\section{Nonabelian examples}\label{sec:examples}

\subsection{The Heisenberg group}

The Heisenberg group $H=\operatorname{UT}_3(\F_3)$ consists of the
$27$ matrices
\[
 (x,y,z)\longleftrightarrow
 \begin{pmatrix}1&x&z\\0&1&y\\0&0&1\end{pmatrix},
 \qquad x,y,z\in\F_3.
\]
All arithmetic below is modulo $3$. Matrix multiplication gives
\begin{equation}\label{eq:group}
 (x,y,z)(u,v,w)=(x+u,y+v,z+w+xv).
\end{equation}
The identity is $(0,0,0)$ and the inverse of $(x,y,z)$ is
$(-x,-y,xy-z)$. With $X=(1,0,0)$ and $Y=(0,1,0)$,
\[
 XY=(1,1,1),\qquad YX=(1,1,0).
\]
Thus the multiplication is noncommutative. Its center, the set of elements commuting with the whole group, consists
of $(0,0,z)$, and the central coordinate is where the multiplication order
is recorded. A commutator measures the failure of two elements to commute.
With the convention $[g,h]=ghg^{-1}h^{-1}$,
\[
 [(x,y,z),(u,v,w)]=(0,0,xv-uy).
\]

\subsection{A homomorphism rule on a group of order 81}\label{sec:example81}

Take $G=H\times(\F_3,+)$ and write elements as $(x,y,z,t)$.
Define
\[
 \alpha(x,y,z,t)=(y,x,xy-z,t),\qquad
 \beta(x,y,z,t)=(0,0,t,t).
\]
For $g=(x,y,z,t)$ and $h=(u,v,w,s)$, the central coordinate of
$\alpha(gh)$ is
\[
 (x+u)(y+v)-(z+w+xv)=xy+uv+uy-z-w.
\]
It equals the central coordinate of $\alpha(g)\alpha(h)$, and the
other coordinates also agree by the displayed formulas. Thus $\alpha$
is an endomorphism. Applying it twice gives $g$, so it is an involutive
automorphism. The map $\beta$ is an endomorphism because its two
nonzero coordinates depend additively on $t$. Its image is central
because $(0,0,t)$ commutes with every element of $H$. Consequently
\begin{equation}\label{eq:hom}
 h(g_1,g_2,g_3)=\alpha(g_1)\beta(g_2)\beta(g_3)
\end{equation}
is a group homomorphism $G^3\to G$: factors from different input
images commute. All three inputs matter, since $\alpha$ is bijective
and $\beta$ is nonconstant.

Prescribe the compatible periodic boundaries
\[
 a(0,n)=(n\bmod3,1,0,1),\qquad
 a(m,0)=\alpha^m(0,1,0,1),
\]
and for $m,n\ge1$ put
\begin{equation}\label{eq:homrec}
 a(m,n)=\alpha(a(m-1,n))\beta(a(m,n-1))\beta(a(m-1,n-1)).
\end{equation}
The predecessors are $(1,0),(0,1),(1,1)$. Corollary~\ref{cor:hom} proves $3$-automaticity.

Writing $a=(x,y,z,t)$ makes the mechanism explicit:
\begin{align*}
 x(m,n)&=y(m-1,n), &y(m,n)&=x(m-1,n),\\
 t(m,n)&=t(m-1,n)+t(m,n-1)+t(m-1,n-1),\\
 z(m,n)&=x(m-1,n)y(m-1,n)-z(m-1,n)\\
       &\hspace{16mm}+t(m,n-1)+t(m-1,n-1).
\end{align*}
The boundary and the first two coordinate equations give
$(x,y)=(n\bmod3,1)$ for even $m$ and $(1,n\bmod3)$ for odd $m$.
In particular, $xy=n\bmod3$ throughout the array.

\subsubsection{The Delannoy coordinate}\label{sec:delannoy}

The \emph{Delannoy number} $D_{m,n}$ counts paths in $\N^2$ from
$(0,0)$ to $(m,n)$ whose steps are $(1,0)$, $(0,1)$ or $(1,1)$.
There is one path along either coordinate axis, so
$D_{m,0}=D_{0,n}=1$. For $m,n\ge1$, partitioning paths according
to their last step gives
\begin{equation}\label{eq:delannoy-definition}
 D_{m,n}=D_{m-1,n}+D_{m,n-1}+D_{m-1,n-1}.
\end{equation}
For example, the three paths to $(1,1)$ are the two possible orders
of the horizontal and vertical steps, and the single diagonal step.
Hence $D_{1,1}=3$.

Prunescu~\cite[Definition~1.2 and Theorem~4.4]{PrunescuCarpets}
studies this recurrence over finite fields with a coefficient on the
diagonal term and proves identities for its residue blocks. Taking
that coefficient to be $1$ gives the Delannoy residues used here.

\begin{figure}[!htbp]
\centering
\begin{minipage}{.475\textwidth}\centering
\textbf{Central coordinate $z$}\\[3pt]
\includegraphics[width=\linewidth]{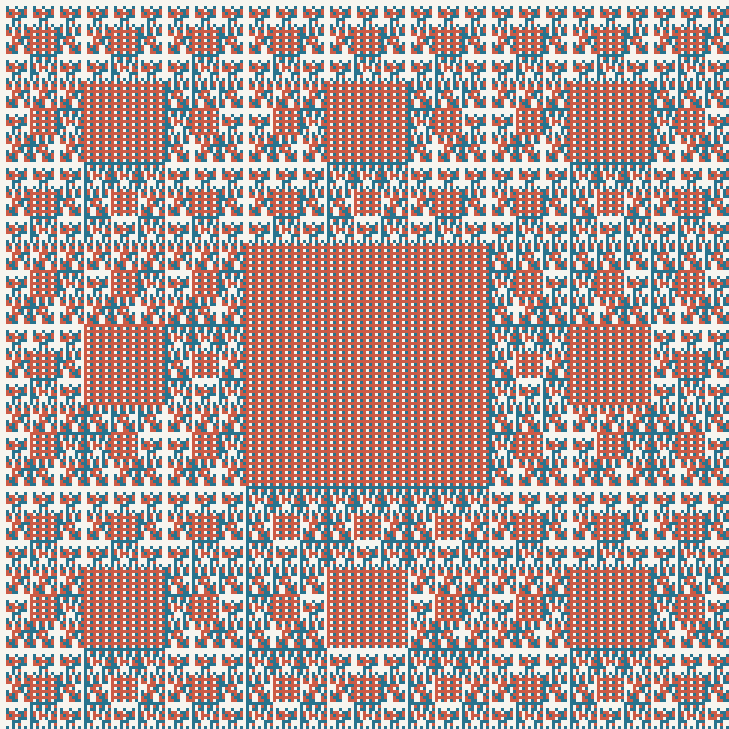}
\end{minipage}\hfill
\begin{minipage}{.475\textwidth}\centering
\textbf{Delannoy coordinate $t$}\\[3pt]
\includegraphics[width=\linewidth]{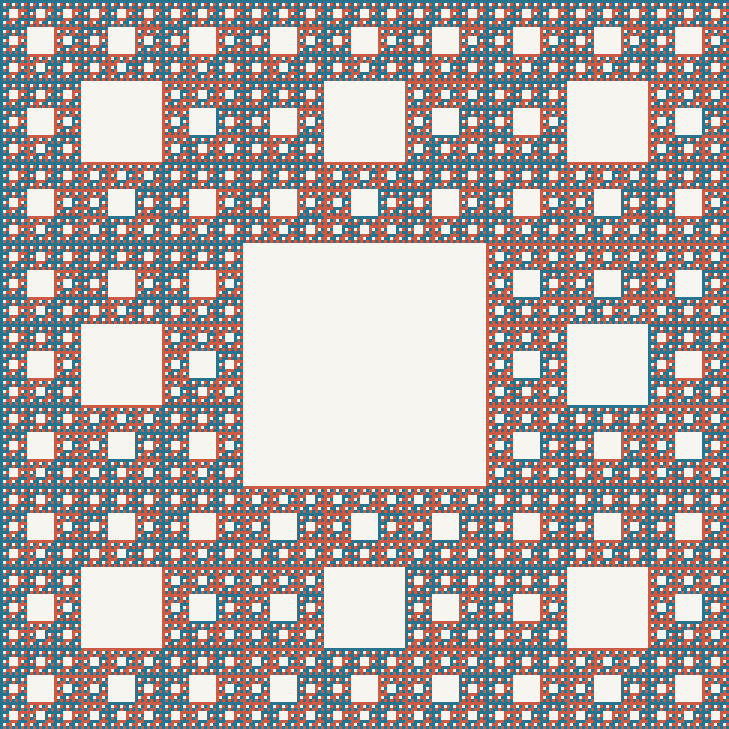}
\end{minipage}
\par\smallskip\legendthree
\caption{The homomorphism example on
$G=\operatorname{UT}_3(\F_3)\times(\F_3,+)$, of order $81$, for
$0\le m,n<243$. The coordinate $n$ increases to the right and $m$
downwards from the upper left corner. Each pixel block shows a residue
modulo $3$. The left panel displays the central coordinate $z$ and the
right panel displays the Delannoy coordinate $t$. Each full group value
is a tuple $(x,y,z,t)$.}\label{fig:hom81}
\end{figure}

\begin{proposition}\label{prop:delannoy-coordinate}
The fourth coordinate of the group recurrence \eqref{eq:homrec} is
$t(m,n)=D_{m,n}\bmod3$.
\end{proposition}
\begin{proof}
Both prescribed boundaries have fourth coordinate $1$.
The maps $\alpha$ and $\beta$ each preserve the input's fourth
coordinate, and multiplication in the last factor of $G$ adds those
coordinates. Thus $t$ satisfies \eqref{eq:delannoy-definition} modulo
$3$ with the same boundary values. Induction on $m+n$ identifies it
with $D_{m,n}\bmod3$.
\end{proof}
We call $t$ the Delannoy coordinate because of this identity. It enters
the equation for $z$ at $(m,n-1)$ and $(m-1,n-1)$. The equation for
$z$ also contains $x(m-1,n)y(m-1,n)-z(m-1,n)$ from the automorphism
$\alpha$. Figure~\ref{fig:hom81} shows these two coordinates over the
square $\{0,\ldots,242\}^2$.

Two interior values witness noncommutativity:
\[
 a(1,1)=(1,1,0,0),\qquad a(1,2)=(1,2,0,2),
 \qquad [a(1,1),a(1,2)]=(0,0,1,0)\ne1.
\]
Thus the solution is not contained in an abelian subgroup.
Among the $59\,049$ values at $0\le m,n<243$, direct evaluation
of the recurrence gives $39$ distinct elements of $G$. The computation
and its checks are described in Section~\ref{sec:verification}.

\clearpage
\subsection{Ordered multiplication with a negative predecessor}\label{sec:example27}

We next combine ordered multiplication in $H$ with a negative
predecessor coordinate. Put $a(0,0)=1$ and prescribe
\[
 a(0,n)=Y\quad(n\ge1),\qquad a(m,0)=X\quad(m\ge1).
\]
For $m,n\ge1$ define
\begin{equation}\label{eq:negativerec}
 a(m,n)=a(m-1,n+1)\,a(m,n-1).
\end{equation}
Its predecessor vectors are $(1,-1)$ and $(0,1)$, both
lexicographically positive. The negative entry is confined to the
second coordinate. The weight $2m+n$ decreases by $1$ for either
predecessor, which gives a direct well-founded construction.

\begin{proposition}
The array in \eqref{eq:negativerec} is $3$-automatic, takes values
outside every abelian subgroup of $H$, and depends on the order of
multiplication.
\end{proposition}
\begin{proof}
Apply Theorem~\ref{thm:groups} with both endomorphisms equal
to the identity. The boundary data are automatic,
and the sign condition holds. Direct computation gives
\[
 a(1,1)=(1,1,0),\qquad a(1,2)=(1,2,0),\qquad
 [a(1,1),a(1,2)]=(0,0,1).
\]
These two interior values cannot lie in a common abelian subgroup.
Already at $(1,1)$ the stated order gives $YX$, while the reversed
order gives $XY$, which are different.
\end{proof}

Multiplication
$H^2\to H$ is not a group homomorphism when $H$ is noncommutative.
This example uses the stronger ordered-product theorem. The order-$81$ example satisfies Corollary~\ref{cor:hom}.

For an elementary view of the proof, write $L=(m-1,n+1)$ and
$D=(m,n-1)$. Its coordinates obey
\[
 x=x_L+x_D,\qquad y=y_L+y_D,\qquad
 z=z_L+z_D+x_Ly_D.
\]
The quotient coordinates $x,y$ are automatic by the abelian theorem
with the allowed negative direction. Their product $x_Ly_D$ is an
automatic affine term. Apply that theorem once more to get automatic
$z$. This application of the abelian theorem supplies the central coordinate
in the group induction.
All boundary central coordinates are zero, but $a(2,2)=(0,2,2)$.
Consequently the nonzero central pattern is produced by the group
multiplication itself. Replacing it by coordinatewise addition would
make the central coordinate identically zero.

\subsubsection{The effect of reversing multiplication}

\begin{figure}[!htbp]
\centering
\begin{minipage}{.475\textwidth}\centering
\textbf{Central coordinate $z$}\\[3pt]
\includegraphics[width=\linewidth]{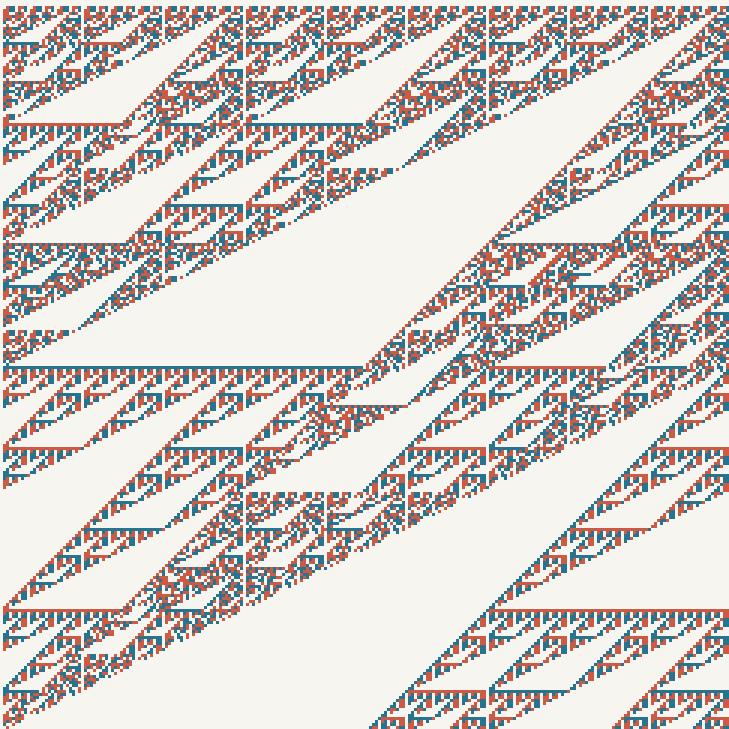}
\end{minipage}\hfill
\begin{minipage}{.475\textwidth}\centering
\textbf{Reversed order: $xy-z$}\\[3pt]
\includegraphics[width=\linewidth]{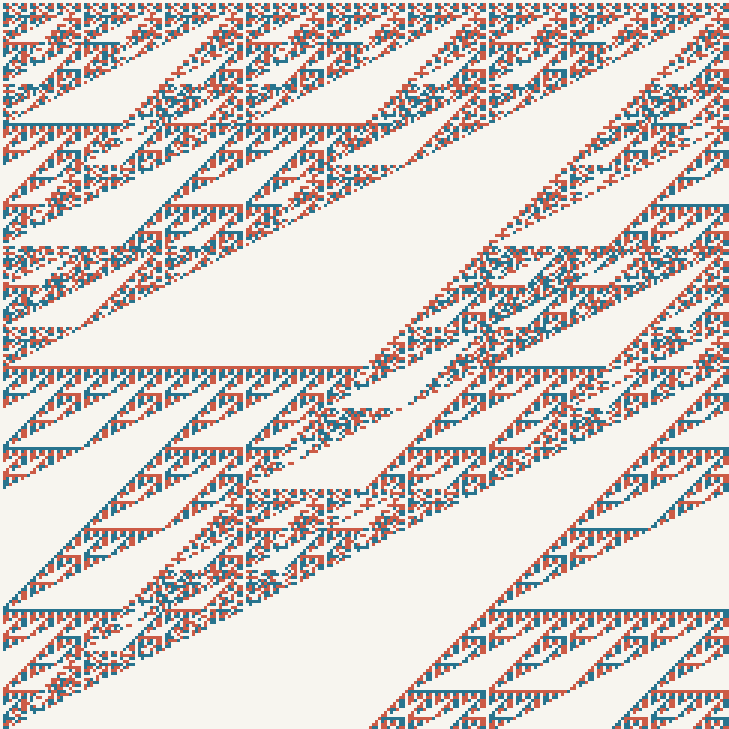}
\end{minipage}
\par\smallskip\legendthree
\caption{Ordered multiplication in $H=\operatorname{UT}_3(\F_3)$,
of order $27$, using
the negative predecessor $(1,-1)$. Coordinates run from $0$ to $242$,
with $n$ to the right and $m$ down. The two panels differ only in the
order of multiplication. Each pixel block shows the residue modulo $3$
of the central coordinate, with $z$ on the left and $xy-z$ on the right.
The other two coordinates agree between the arrays.}\label{fig:ordered27}
\end{figure}

Let $a^{\rm rev}$ have the same boundary and reversed product in
\eqref{eq:negativerec}. Define
\[
 \tau(x,y,z)=(x,y,xy-z).
\]
Calculation gives $\tau(gh)=\tau(h)\tau(g)$ and $\tau(X)=X$,
$\tau(Y)=Y$. Induction on $2m+n$ proves the exact identity
\[
 a^{\rm rev}(m,n)=\tau(a(m,n)),\qquad
 z^{\rm rev}(m,n)=x(m,n)y(m,n)-z(m,n).
\]
The two arrays have the same coordinates $x,y$. Among the
$59\,049$ positions in $0\le m,n<243$, direct evaluation finds
$17\,737$ positions where their central coordinates differ. All $27$
elements of $H$ occur in this square. Figure~\ref{fig:ordered27}
compares the two central coordinates.

\clearpage
\section{Exact computations on the displayed squares}\label{sec:verification}

The figures and counts concern the square
$S=\{(m,n):0\le m,n<243\}$, which contains $59\,049$ positions.
The computation evaluates the specified recurrences at all positions
in this square and records group elements as residue tuples. Distinct
values are counted as distinct tuples, without identifying values by
an automorphism or any other equivalence relation.
Table~\ref{tab:finite-counts} records the counts.

\begin{table}[!htbp]
\centering
\begin{tabular}{llr}
\toprule
Array & Quantity on $S$ & Count\\
\midrule
Order-$81$ example & Distinct values in $G$ & $39$\\
Order-$27$ example & Distinct values in $H$ & $27$\\
Order-$27$ example & Positions changed by reversing multiplication & $17\,737$\\
\bottomrule
\end{tabular}
\caption{Exact counts on $S=\{0,\ldots,242\}^2$, containing
$59\,049$ positions. Distinct values are counted as group elements.
The final row compares the two multiplication orders at each position.}
\label{tab:finite-counts}
\end{table}

The last count compares the full group values. Since reversing the
product preserves $x$ and $y$, it is also the count of positions where
the central coordinate changes.

For the order-$81$ example, evaluation proceeds by increasing $m$ and,
within a row, increasing $n$. The three predecessors in
\eqref{eq:homrec} have already been evaluated. For the negative
predecessor in \eqref{eq:negativerec}, a rectangular truncation would
omit required values at $n+1$. Instead, for a square of side $N$,
row $m$ is evaluated through $n=2N-2-m$, for $0\le m<N$.
The preceding row is one entry longer, so its value at $n+1$ is
available. The current row's value at $n-1$ is also available.
Because $2N-2-m\ge N-1$, the computation covers the entire requested
square before cropping. This proves that the computed entries
agree with the infinite recurrence at all $N^2$ positions used for
each count.

The repository cited in footnote~\ref{fn:repository} includes the script
\path{examples/nonabelian_examples.py}, which performs these evaluations
using exact arithmetic modulo $3$. It writes the images,
small arrays and counts under \path{examples/generated/nonabelian/},
including \path{verification.json}.
It also checks the two endomorphism identities and centrality of the
image of $\beta$ on all $81^2=6\,561$ ordered pairs in each case.
Associativity in $H$ is checked on all $27^3=19\,683$ ordered triples.
An independent computation ordered by the weight $2m+n$ checks
the negative recurrence on $0\le m,n<35$. The identity with $\tau$
is checked at all $59\,049$ displayed positions.

We also check the Delannoy coordinate independently. Define
\[
 F(X,Y)=\sum_{m,n\ge0}D_{m,n}X^mY^n,\qquad Q=1-X-Y-XY.
\]
The recurrence cancels the coefficients of $QF$ off the axes.
On a positive axis, the unit boundary gives $1-1=0$, and the
constant coefficient is $1$. Therefore $F=1/Q$.
Over $\F_3$, $Q^3=Q(X^3,Y^3)$, so
\[
 F(X,Y)=Q(X,Y)^2 F(X^3,Y^3),\qquad
 Q^2=\sum_{0\le i,j<3} B_{ij}X^iY^j,
 \quad B=\begin{pmatrix}1&1&1\\1&0&2\\1&2&1\end{pmatrix}.
\]
Coefficient extraction gives
$D_{3m+i,3n+j}\equiv B_{ij}D_{m,n}\pmod3$.
Iteration computes a residue as the product of entries $B_{ij}$
selected by successive ternary digit pairs. The script compares this
formula with the coordinate $t$ at every displayed position.

These computations establish the counts on $S$ and check the
implementation of the examples. Automaticity of the infinite arrays follows from the
theorems and proofs in Sections~\ref{sec:negative}--\ref{sec:examples}.

\section{An exact Delannoy count modulo 27}\label{sec:mod27-counts}

The Delannoy numbers defined in Section~\ref{sec:delannoy} also permit
exact counts of nonzero residues. For $r\ge0$, put
\[
 N_{27}(r)=|\{(m,n):0\le m,n<3^r,\ D_{m,n}\not\equiv0\pmod{27}\}|.
\]
Laurberg \cite[Eq.~(10.47), p.~164]{Laurberg} proposes the following
formula from computed values. We give a certificate consisting
of a transition table and identities between polynomial and integer
vectors. These identities prove the formula for all $r\ge4$.

\begin{theorem}\label{thm:mod27-count}
For every $r\ge4$,
\[
 N_{27}(r)=8^{r-4}(18r^2+471r+3993).
\]
Consequently, the number of zero entries in that square is
$9^r-N_{27}(r)$.
\end{theorem}
\begin{proof}
Work over $R=\Z/27\Z$. The generating function of the Delannoy
array is $1/Q$, where $Q=1-X-Y-XY$, as proved in
Section~\ref{sec:verification}. The construction in
Theorem~\ref{thm:division} uses states $P/Q^9$, starts with
$P=Q^8$, and has transitions
\[
 P\longmapsto\Lambda_{i,j}(P Q^{18}),\qquad 0\le i,j<3.
\]
Enumerating the reachable polynomials gives $1362$ states.
The file \path{research/generated/mod27_rational_machine.json}
in the repository cited in footnote~\ref{fn:repository} records their
coefficient arrays.
The script \path{research/verify_mod27.py} verifies the initial
polynomial, reachability, constant-coefficient outputs, and all
$9\cdot1362=12\,258$ transition identities over $R[X,Y]$.
Each polynomial has degree at most $9$ in each variable, so the
state bound of Theorem~\ref{thm:division} proves termination of
this enumeration.

Index the states starting from $q_0=0$, and define the integer matrix
\[
 M_{q,q'}=|\{(i,j):\delta(q,(i,j))=q'\}|.
\]
Let $e_0$ be the column vector indexed by the $1362$ states, with
entry $1$ at state $0$ and $0$ elsewhere. Let $v$ be the column
vector with $v_q=1$ when state $q$ outputs a nonzero residue
modulo $27$, and $v_q=0$ otherwise. Each coordinate pair below $3^r$ has one padded
digit representation of length $r$, so
\[
 N_{27}(r)=e_0^{\mathsf T}M^r v.
\]
The certificate stores the eight integer vectors $v,Mv,\ldots,M^7v$
and checks, in all $1362$ coordinates, the identity
\begin{equation}\label{eq:mod27-certificate}
 M^7v-24M^6v+192M^5v-512M^4v=0.
\end{equation}
Multiplying by $e_0^{\mathsf T}M^{r-4}$ gives, for $r\ge4$,
\[
 N_{27}(r+3)-24N_{27}(r+2)+192N_{27}(r+1)-512N_{27}(r)=0.
\]
The same certificate gives
\[
 N_{27}(4)=6165,\qquad N_{27}(5)=54384,\qquad N_{27}(6)=477888.
\]
The proposed formula has these three values. It satisfies the displayed
recurrence because applying the operator $f(r)\mapsto f(r+1)-8f(r)$
three times to $8^r$ times a quadratic polynomial gives zero.
The recurrence and initial values therefore prove the formula.
\end{proof}

As an independent check, the script also evaluates the Delannoy recurrence
directly on a $243\times243$ square and compares all $59\,049$
residues with the automaton.

\begin{corollary}\label{cor:mod27-density}
The proportion of entries divisible by $27$ in
$\{0,\ldots,N-1\}^2$ tends to $1$ as $N\to\infty$.
\end{corollary}
\begin{proof}
For $N=3^r$, the nonzero proportion is
$8^{-4}(18r^2+471r+3993)(8/9)^r$, which tends to zero.
For $3^r\le N<3^{r+1}$, the number of nonzero entries is at most
$N_{27}(r+1)$ and the square has at least $9^r$ entries. Thus its
nonzero proportion is at most $N_{27}(r+1)/9^r$, which also tends
to zero.
\end{proof}

\subsection{State counts for the first three powers of 3}

For the numerical Delannoy arrays, Table~\ref{tab:delannoy-states}
compares certified minimal state counts with a bound from the rational
case of Theorem~\ref{thm:division}. All machines in this comparison read
least significant ternary digits first and output the full residue.
To count polynomial states, put $E=3^{e-1}$. The numerator series $1$
has kernel $\{1,0\}$. In the construction of Theorem~\ref{thm:division},
the polynomial multiplying the zero series can be discarded.
The remaining state is a single polynomial of degree at most $E$ in
each variable, with $(E+1)^2$ coefficients in $\Z/3^e\Z$.
This gives the bound $3^{e(E+1)^2}$.

\begin{table}[!htbp]
\centering
\begin{tabular}{rrr}
\toprule
Modulus & Minimal states & Construction bound\\
\midrule
$3$ & $3$ & $3^4$\\
$9$ & $58$ & $3^{32}$\\
$27$ & $1362$ & $3^{300}$\\
\bottomrule
\end{tabular}
\caption{Exact minimal state counts for the two-dimensional numerical
Delannoy array modulo $3^e$, for $e=1,2,3$, using least-significant-digit
automata with ternary digit pairs and zero padding. The last column is
a bound on the polynomial states before reachability and minimization.}
\label{tab:delannoy-states}
\end{table}

Reachability and refinement of states by their outputs and successor
classes certify the minimal counts. These three cases do not determine
the minimal state count as a function of $e$.

\section{Further questions}\label{sec:questions}

Theorem~\ref{thm:mixed} guarantees automaticity when at most one
coordinate contains predecessor entries of both signs. The fixed
obstruction shows that two such coordinates can already produce a
nonautomatic array. It remains to determine which predecessor sets
with two or more such coordinates guarantee automaticity for every
permitted rule and every compatible automatic datum.

Theorem~\ref{thm:asymmetric} gives an exact criterion for one specified
rule and seed. Its positive cases do not establish automaticity for
all rules on the same predecessor set. Understanding which features
of a boundary equation produce a multiplicative substitution such as
$X\mapsto X^{bc}$ could give further obstruction criteria.

The verified search of Theorem~\ref{thm:verified-search} gives a
terminating construction under the geometric theorems. The proof gives
no bound on its number of computational steps in terms of the input
rule and automata. A direct construction from the matrix division proof
could give bounds in terms of the group, predecessors and input automata.
Such bounds could be compared, where their hypotheses apply, with
the explicit bounds of \cite[Theorems~1--3]{RY26} for algebraic
series and diagonals of rational functions modulo prime powers.

Theorem~\ref{thm:mod27-count} gives an exact count of nonzero Delannoy
residues modulo $27$. For the family modulo $3^e$, with $e\ge4$, the
results here do not determine how the counting recurrences and their
polynomial factors depend on $e$, or whether the corresponding minimal
automata admit a uniform description. The three state counts in
Table~\ref{tab:delannoy-states} alone do not determine that dependence.
Laurberg \cite[Section~10.3, Eq.~(10.49)]{Laurberg} conjectures that,
for every prime $p$ and every $e\ge1$, the number of nonzero residues
modulo $p^e$ in the square of side $p^r$ is eventually the
corresponding count modulo $p$ times a polynomial in $r$ of degree
$e-1$. The exact-degree assertion requires an exception at $p=2$
when $e>1$. Indeed, the unit boundary and the Delannoy recurrence
show by induction that every $D_{m,n}$ is odd. The count modulo
$2^e$ is therefore $4^r$ for every $e\ge1$, so the multiplier
relative to the count modulo $2$ is the constant polynomial $1$.
For $p=3$ and $e=3$, Theorem~\ref{thm:mod27-count} establishes
the eventual formula predicted by Laurberg, with count $8^r$
modulo $3$ and a quadratic multiplier for $r\ge4$.

\section*{Acknowledgments}

This paper was produced in an AI-assisted workflow. AI assistants helped explore candidate constructions, draft proofs, and write software for numerical experiments.

\end{document}